\documentclass[11pt,authoryear]{elsarticle}
\usepackage[a4paper,margin=2cm]{geometry}
\usepackage{graphicx}
\usepackage{amsmath,amsthm,amssymb}
\usepackage{booktabs}
\usepackage{threeparttable}
\usepackage{tabularx}
\usepackage{array}
\usepackage{multirow}
\usepackage{makecell}
\usepackage{siunitx}
\usepackage{setspace}
\usepackage{microtype}
\usepackage{xcolor}
\usepackage{longtable}
\usepackage{enumitem}

\usepackage{hyperref}
\hypersetup{
    hidelinks,
    colorlinks=true,
    allcolors=black,
    pdfstartview=Fit,
    breaklinks=true
}
\usepackage[nameinlink,noabbrev]{cleveref}
\newcommand{\sig}{\rlap{\textsuperscript{*}}}
\usepackage{placeins}

\begin{document}
\begin{frontmatter}

\title{Automated Priority Rule Design for the Resource-Constrained Project Scheduling Problem: A Large Language Model-Guided Population-Based Search}

\author[mymainaddress]{Jingyu Luo}
\ead{jingyu.luo@mcgill.ca}

\author[mysecondaryaddress,mythirdaddress,myfourthaddress]{Mario Vanhoucke \corref{mycorrespondingauthor}}
\cortext[mycorrespondingauthor]{Corresponding author}
\ead{mario.vanhoucke@ugent.be}
		
\author[myfifthaddress,mysixthaddress]{Jos\'e Coelho}
\ead{jose.coelho@uab.pt}

\address[mymainaddress]{COSMO – Stochastic Mine Planning Laboratory, McGill University, 3450 Rue University, H3A 0E8 Montreal (Canada)}
\address[mysecondaryaddress]{Ghent University, Tweekerkenstraat 2, 9000 Ghent (Belgium)}
\address[mythirdaddress]{University College London, 1 Canada Square, E14 5AA London (United Kingdom)}
\address[myfourthaddress]{Vlerick Business School, Reep 1, 9000 Ghent (Belgium)}
\address[myfifthaddress]{Universidade Aberta, Rua da Escola Polit\'ecnica, 147, 1269-001 Lisbon (Portugal)}
\address[mysixthaddress]{INESC TEC, Campus da FEUP, Rua Dr.Roberto Frias, 4200 - 465 Porto (Portugal)}

\begin{abstract}
The objective of the resource-constrained project scheduling problem (RCPSP) is to minimize makespan while satisfying precedence and renewable-resource constraints. Priority-rule heuristics construct feasible schedules with low computational cost and explicit decision logic, making them widely used in practice and an attractive alternative to more computationally intensive methods. However, no traditional rule performs consistently well across projects, and researchers have therefore investigated automated priority-rule design. Genetic programming (GP) hyper-heuristics have been the predominant approach to this task, but evolving a high-performing rule may require evaluating many candidate rules on the training projects. Recent advances in large language models (LLMs) make it possible to generate and iteratively revise priority rules using performance feedback.

This paper presents an LLM-guided population-based framework for automated priority-rule design. During the offline search, an LLM generates and revises candidate rules, while schedule quality on training projects determines candidate fitness and guides subsequent revisions. At the end of the search, the best rule is returned and applied directly to unseen projects without further search or LLM calls.

Experiments show that the LLM-designed rules outperform traditional single rules across the test sets and outperform GP-designed rules obtained under comparable search effort, while remaining competitive with rules obtained from a substantially larger GP search. On large projects, selected LLM-designed rules outperform all considered traditional rules and two genetic algorithm configurations on highly parallel projects. An ablation study examines how changes to the main search components affect performance, while rule analyses describe the structure and decision behavior of the LLM-designed rules.
\end{abstract}


\end{frontmatter}

\section{Introduction}
\label{sec:introduction}

The resource-constrained project scheduling problem (RCPSP) concerns the nonpreemptive scheduling of project activities subject to precedence relations and limited renewable resources, with the objective of minimizing the project makespan. Researchers have addressed the problem using exact methods, metaheuristics, and priority-rule heuristics. Exact and metaheuristic methods can obtain high-quality schedules, but their computational requirements may limit their use for large projects or in settings that require rapid scheduling and rescheduling. Priority-rule heuristics provide a simpler alternative by using explicit ranking rules to construct feasible schedules at low computational cost. Their ease of implementation and transparent decision logic make them particularly attractive when projects are large or schedules must be generated repeatedly under tight time constraints \citep{djumic2021ensembles,luo2022efficient}. 

The performance of a priority rule varies across project instances because different rules emphasize different timing, precedence, and resource information. Consequently, no traditional rule performs consistently well across diverse projects. Designing new rules manually requires substantial domain knowledge and repeated experimentation, and the resulting rules may still perform poorly on projects that differ from those considered during design. These limitations have motivated methods that automatically select existing priority rules or generate new ones.

For the standard RCPSP, automated priority-rule methods have mainly followed two routes: genetic programming hyper-heuristics (GPHH, \citet{koza1992genetic}) for direct rule generation and supervised learning from labeled scheduling data. GPHH has been the principal approach to generating new priority rules. In these methods, genetic programming (GP) evolves explicit scoring expressions according to the quality of the schedules they construct on training projects \citep{chand2018use,djumic2018evolving,luo2022efficient}. Supervised approaches instead use labeled scheduling data either to select existing rules or to learn activity-ranking models \citep{guo2021automatic,luo2026automated}. GPHH may require evaluating many candidate rules before a high-performing rule is found, while supervised methods depend on the availability and quality of labeled data. The transparency of supervised models also depends on the model class.

Large language models (LLMs) offer another route to automatic heuristic design because they can generate executable priority functions that use specified scheduling inputs and return a priority score. Evaluation-guided LLM search has shown promising results in combinatorial optimization and production scheduling \citep{romera2024mathematical,liu2024evolution,ye2024reevo,yu2026automated}. The recent fifty-year review of RCPSP research by \citet{artigues2026fifty} identifies machine-learning-based heuristic design as an emerging direction but does not discuss LLM-based priority-rule design. Within project scheduling, \citet{kim2026vlcea} use LLM-generated, instance-specific priority scores to guide evolutionary schedule search, while \citet{tian2026guiding} use GP-derived knowledge to guide an LLM for online activity--mode decisions. These studies show that LLMs can support project-scheduling decisions, but their use for automatic priority-rule design for the standard RCPSP remains limited. In both studies, the LLM contributes to solving individual project instances. By contrast, the present study targets an executable priority function designed offline and reused across unseen projects. For this design setting, evidence remains limited on how to organize the search and how its performance compares with established automatic rule-design methods under controlled search effort. The scheduling logic encoded by the resulting rules has also received little attention.

This study proposes an LLM-guided population-based framework for offline priority-rule design for the standard RCPSP. All candidate rules follow the same executable function format and use the same predefined scheduling inputs. The LLM first generates an initial population of candidate rules. Their fitness is determined by the schedules they produce on the training projects, and the resulting feedback guides structured rule revisions in later generations. Selected rules are carried forward, and reference rules are drawn from a Hall of Fame that balances performance and diversity. Each run returns one deterministic priority rule that can be inspected and applied directly to unseen projects. Deployment requires no further rule-design search or LLM calls.

This study makes three contributions. First, it develops a population-based LLM search framework tailored to the offline design of reusable and executable RCPSP priority rules. Second, it provides systematic computational evidence on the performance of the designed rules under different training settings and with different LLMs, as well as their generalization across benchmark families and project sizes. The comparisons include traditional and automatically designed priority rules, as well as genetic algorithm results for large projects. Third, it uses an ablation study to examine the contribution of the search components to performance and analyzes the structure and decision behavior of the generated rules.

The remainder of the paper is organized as follows. \Cref{sec:background} introduces the RCPSP and reviews automatic priority-rule design and LLM-guided heuristic design. \Cref{sec:framework} presents the proposed framework, and \Cref{sec:experimental_design} describes the benchmark data and experimental protocol. \Cref{sec:results} reports the computational results, while \Cref{sec:rule_analysis} analyzes the generated rules and discusses the resulting priority-rule design insights. Finally, \Cref{sec:conclusion} concludes the study and outlines directions for future research.

\section{Problem background and related work}
\label{sec:background}

\subsection{The resource-constrained project scheduling problem}
\label{sec:problem_description}

The RCPSP considers the nonpreemptive scheduling of project activities subject to precedence relations and limited renewable resources \citep{slowinski1981multiobjective}. Let \(N\) denote the number of non-dummy activities, let \(\mathcal{A}=\{0,1,\ldots,N,N+1\}\) denote the complete activity set, and let \(R\) denote the set of renewable resources. Activities 0 and \(N+1\) are dummy start and finish activities, respectively, with zero duration and zero resource requirements. Each non-dummy activity \(i\in\{1,\ldots,N\}\) has a fixed duration \(d_i>0\) and requires \(r_{i,k}\geq 0\) units of each renewable resource \(k\in R\) throughout its execution. Each resource \(k\) has a constant availability \(a_k\).

The precedence relations form a directed acyclic activity-on-node network. Let \(IP_i\) and \(S_i\) denote the sets of immediate predecessors and immediate successors of activity \(i\), respectively. A precedence relation from activity \(i\) to activity \(j\) requires \(j\) to start only after \(i\) has finished. A schedule assigns a start time \(s_i\) to each activity \(i\in\mathcal{A}\). It is feasible if all precedence relations are respected and the total demand for each renewable resource does not exceed its availability at any time.

Using this notation, the RCPSP can be formulated as
\begin{align}
\min\quad & C_{\max}=s_{N+1}, \label{eq:rcpsp_objective}\\
\text{s.t.}\quad & s_i+d_i\leq s_j, && \forall i\in\mathcal{A},\ \forall j\in S_i, \label{eq:precedence_constraint}\\
& \sum_{\substack{i\in\mathcal{A}:\\s_i\leq \tau<s_i+d_i}} r_{i,k}\leq a_k,
&& \forall k\in R,\ \forall \tau\geq 0,
\label{eq:resource_constraint}\\
& s_0=0,\qquad s_i\geq 0, && \forall i\in\mathcal{A}. \label{eq:start_time_constraint}
\end{align}

The objective in \Cref{eq:rcpsp_objective} minimizes the project makespan, which equals the start time of the zero-duration finish activity. Constraints~\eqref{eq:precedence_constraint} enforce the precedence relations, while Constraints~\eqref{eq:resource_constraint} ensure that the concurrent demand for each renewable resource does not exceed its availability. Constraint~\eqref{eq:start_time_constraint} fixes the dummy start activity at time zero and imposes nonnegative start times.

\subsection{Priority rules and schedule generation}
\label{sec:priority_rules_sgs}

Priority-rule heuristics combine a priority rule with a schedule generation scheme (SGS). The SGS constructs a feasible schedule by successively extending a partial schedule. At decision point \(t\), let \(E_t\) denote the set of eligible non-dummy activities. The precise eligibility conditions and the procedure for assigning feasible start times depend on the SGS, while the priority rule determines which activity is selected from \(E_t\). A priority rule assigns a numerical priority value to each eligible activity using information available at the current decision point. This information may include activity attributes, such as activity duration, temporal position, precedence-network position, and resource requirements, as well as project-level or schedule-construction information. We refer to the variables through which this information is supplied to the priority rule as \textit{rule inputs}. The rule also specifies whether the activity with the smallest or largest priority value is selected. A tie-breaking rule is applied when multiple activities receive the same priority value \citep{kolisch1996serial}.

The two principal schedule generation schemes for the RCPSP are the serial schedule generation scheme (SSGS) and the parallel schedule generation scheme (PSGS). Both start from an initial partial schedule and extend it while maintaining precedence and renewable-resource feasibility. Their main distinction is that SSGS proceeds by activity incrementation, whereas PSGS proceeds by time incrementation \citep{kolisch1996serial,kolisch2014shifts}.

Under SSGS, one non-dummy activity is selected at each stage. The eligible set contains the unscheduled activities whose immediate predecessors have already been scheduled. After the priority rule selects one activity, the activity is assigned its earliest precedence- and resource-feasible start time. The partial schedule is then updated, and the procedure continues until all non-dummy activities have been scheduled.

Under PSGS, schedule construction proceeds through successive scheduling times. At a given scheduling time, an activity is eligible if all of its immediate predecessors have been completed and its resource requirements can be accommodated by the currently available resources. Activities are selected successively from the eligible set and started at the current scheduling time until no additional eligible activity can be started. The scheduling time then advances to the earliest completion time among the activities in progress, after which resource capacity is released, and the eligible set is updated. This process continues until all non-dummy activities have been scheduled.

The difference between the two schemes can therefore be stated operationally. Each SSGS selection inserts one activity into the partial schedule at its earliest feasible start time. PSGS instead holds the scheduling time fixed while starting eligible activities and advances time only when no further activity can be started. Because the two schemes generate different eligible sets and update the partial schedule in different ways, the same priority rule may select different activities and ultimately produce different schedules under SSGS and PSGS.

\subsection{Automatic priority-rule design for the RCPSP}
\label{sec:automatic_rule_design}

Automated priority-rule methods for the RCPSP have mainly followed two routes: direct rule generation through GPHH and supervised learning from labeled scheduling data. GPHH belongs to the broader hyper-heuristic literature, where the search operates on heuristics rather than directly on problem solutions \citep{burke2013hyper}. In the RCPSP setting, GP searches over priority-rule expressions rather than project schedules. Supervised methods use labels either to select an existing rule for a project or to learn activity-ranking information from previously generated schedules.

In GPHH, a candidate priority rule is generally represented as an expression tree. Its terminal nodes contain activity attributes and numerical constants, while its internal nodes contain mathematical operators. Each candidate is embedded in an SGS and applied to the training projects. The quality of the resulting schedules determines its fitness, and selection, crossover, and mutation are used to evolve the rule population. Early work includes \citet{frankola2008evolutionary}, followed by more extensive studies by \citet{djumic2018evolving}, \citet{chand2018use}, and \citet{luo2022efficient}. These studies show that rules evolved on training projects can retain effective performance on unseen test projects and outperform traditional single priority rules across various RCPSP benchmark datasets.

Later studies have extended GPHH in several ways. Rollout justification examines additional scheduling choices before a scheduling decision is finalized and can improve the schedule quality achieved by evolved rules \citep{chand2019evolving,dhumic2022using}. Ensemble approaches combine the decisions of multiple GP-designed rules rather than deploying only one rule \citep{djumic2021ensembles}. \citet{chen2023guided} introduced attribute-node activation encoding and guided genetic operations based on estimated contributions of attribute nodes and subtrees.

Other studies have focused on reducing the computational requirements of GPHH. \citet{luo2022efficient} examined training-set selection, duplicate removal, and fitness-evaluation procedures. The results showed that carefully selected training instances and duplicate removal can improve search efficiency while retaining competitive test-set performance. \citet{luo2023automated} subsequently introduced surrogate models to approximate candidate fitness, accelerate fitness evaluation, and shorten the overall training process. Nevertheless, GPHH can still require substantial computational effort before a high-performing rule is found, because many candidate rules must be assessed and each exact fitness evaluation involves constructing schedules across the training set.

Supervised learning has been used both to select a rule for an entire project and to learn activity rankings. \citet{guo2021automatic} used decision trees to select a suitable traditional priority rule from a predefined collection based on project-level indicators. This approach selects one rule for each project but does not generate a new priority function.

Other supervised models learn activity-ranking information from high-quality schedules. \citet{teichteil2023fast} trained a graph neural network regression model to predict activity start times using schedules obtained from a constraint-programming solver. Activities are ordered according to the predicted start times, after which SSGS constructs a feasible schedule. The predictions therefore induce an activity ordering, although they do not form an explicit priority function expressed through a compact set of rule inputs.

More recently, \citet{luo2026automated} investigated regression-based heuristic design for the RCPSP. Several regression algorithms and ensemble models were trained to learn relationships between activity and project information and priority values derived from high-quality schedules. Unlike the start-time predictions described above, the fitted models are used directly as scoring rules within SSGS and PSGS to rank activities in unseen projects. Across various RCPSP benchmark datasets, the resulting heuristics outperform traditional single priority rules and achieve competitive or superior performance relative to GP-designed rules on selected datasets.

The two routes have complementary strengths and limitations. GPHH does not require labeled schedules and evaluates candidate rules directly according to the quality of the schedules they construct. Its repeated fitness evaluations, however, can require substantial computational effort. Supervised methods avoid population-based rule evolution once suitable labeled data are available, but their training depends on high-quality schedules or rule-performance labels. The transparency of the learned output also depends on the model class: symbolic regression and individual decision trees can provide explicit rules, whereas neural networks and some ensemble models are more difficult to inspect or express as compact priority functions. These trade-offs motivate methods that use schedule quality directly while retaining an explicit and reusable priority-rule representation.

\subsection{Large language model-guided heuristic design and research gap}
\label{sec:llm_heuristic_design}

Recent advances in the code-generation capabilities of large language models (LLMs) have broadened the scope of automatic heuristic design. GP-based methods generally search over expressions assembled from predefined terminals and functions, whereas LLM-based methods can generate complete executable heuristics from natural-language task specifications. The prompt, required interface, and validation procedure constrain the admissible outputs. Even so, LLM-based methods can search over executable program structures without defining the entire search space in advance through a fixed expression grammar. FunSearch combines LLM-based program generation with external performance evaluation and evolutionary selection \citep{romera2024mathematical}. Evolution of Heuristics (EoH) evolves both natural-language descriptions of heuristic ideas and their executable implementations \citep{liu2024evolution}. ReEvo, a reflective-evolution framework, uses performance feedback to guide revisions of heuristic programs \citep{ye2024reevo}. In these methods, candidate heuristics are evaluated on training problems, and the resulting performance feedback is incorporated into later rounds of LLM-based heuristic generation and revision.

The same evaluation-guided pattern has increasingly been adopted in production scheduling. \citet{yu2026automated} proposed a framework that combines LLM-based heuristic generation, performance evaluation, and evolutionary search for flow shop, job shop, and open shop scheduling. For dynamic fuzzy job shop scheduling, \citet{huang2026automatic} used population-level information to support the iterative generation and improvement of dispatching rules. \citet{huang2026automated} subsequently introduced evolutionary textual gradients, which aggregate feedback across a population of candidate dispatching rules and convert the identified improvement directions into instructions for later revisions. DGEvo uses both performance information and rule-structure information to guide LLM-based heuristic evolution for a job shop scheduling problem with fixed preventive maintenance \citep{liu2026dgevo}. These studies show that LLMs can generate and improve executable dispatching and sequencing rules across a range of production-scheduling problems.

LLMs have also been integrated into broader algorithms for scheduling problems. \citet{li2026memetic} incorporated LLM-generated search strategies into a memetic algorithm for multi-objective flexible job shop scheduling with variable processing speeds. \citet{wang2025multiagent} use multiple LLM agents to formulate scheduling problems, initialize solution populations, and carry out evolutionary improvement and evaluation. The framework is tested on several scheduling problems. \citet{wang2026large} combined LLM assistance with an evolutionary order-dispatching method for multiple agile earth-observing satellite scheduling. \citet{zhao2026dga2d} proposed DGA\(_2\)D, in which LLMs generate and revise code implementations for functional operators and modify how these operators are connected. A directed path through the resulting graph specifies a complete candidate algorithm, and the RCPSP is included among the problems used to evaluate the framework. In these studies, LLMs contribute search strategies, dispatching logic, algorithm components, or search coordination within a larger optimization procedure.

Within project scheduling, recent studies have applied LLMs in two distinct ways. For the standard RCPSP, \citet{kim2026vlcea} proposed a method in which an LLM generates validated activity-priority scores separately for each project instance, and these scores guide a permutation-based evolutionary algorithm. Tests on projects of different sizes, up to 120 activities, show that the generated guidance can improve the evolutionary search, although its effectiveness depends on how the guidance is incorporated into the algorithm. For dynamic multi-mode project scheduling, \citet{tian2026guiding} extract knowledge from GP-designed priority rules and use it to guide an LLM that selects activity--mode pairs at online decision points. Several forms of this guidance improve scheduling performance relative to the same LLM without such guidance, with different trade-offs between solution quality and token consumption. The two studies therefore assign the LLM different roles: \citet{tian2026guiding} transfer GP-derived knowledge to an LLM that supports online activity--mode decisions, whereas \citet{kim2026vlcea} use LLM-generated, instance-specific scores to guide an evolutionary schedule search. In both cases, the LLM contributes to solving a particular project instance.

These studies provide evidence that LLMs can support online decisions and instance-specific search in project scheduling. They do not, however, examine offline search for a fixed priority rule that can be reused across unseen projects. This different design target shifts the computational effort to an offline rule-design stage. Once selected, the rule can be inspected and applied to new projects without a new rule-design search or additional LLM calls.

The resulting rule retains the main practical advantages of priority-rule heuristics, namely simple implementation and low computational cost \citep{djumic2021ensembles}. These properties are particularly valuable for large projects, where repeated per-instance optimization may be costly. The project-scheduling studies reviewed above do not examine reusable-rule deployment on large projects. Evidence across benchmark families and different LLMs is also limited, and little attention has been paid to the inputs, structures, and scheduling decisions encoded by the generated rules.

To address these gaps, this study proposes an LLM-guided population-based framework for offline priority-rule design for the standard RCPSP. Each search run returns an explicit priority function that can be applied to unseen projects without further LLM calls. To assess the framework's performance and robustness, the experiments vary the training data and the LLM used for rule generation and test the resulting rules across benchmark families and project sizes. We compare these rules with traditional, GP-designed, and regression-based priority rules and include genetic algorithm results for the large-project datasets. An ablation study examines the contribution of the search components to performance, while the rule analyses show how the generated rules are structured and how they behave during schedule construction.

\section{LLM-guided priority-rule design framework}
\label{sec:framework}

\begin{figure}[htbp]
\centering
\includegraphics[width=0.9\textwidth]{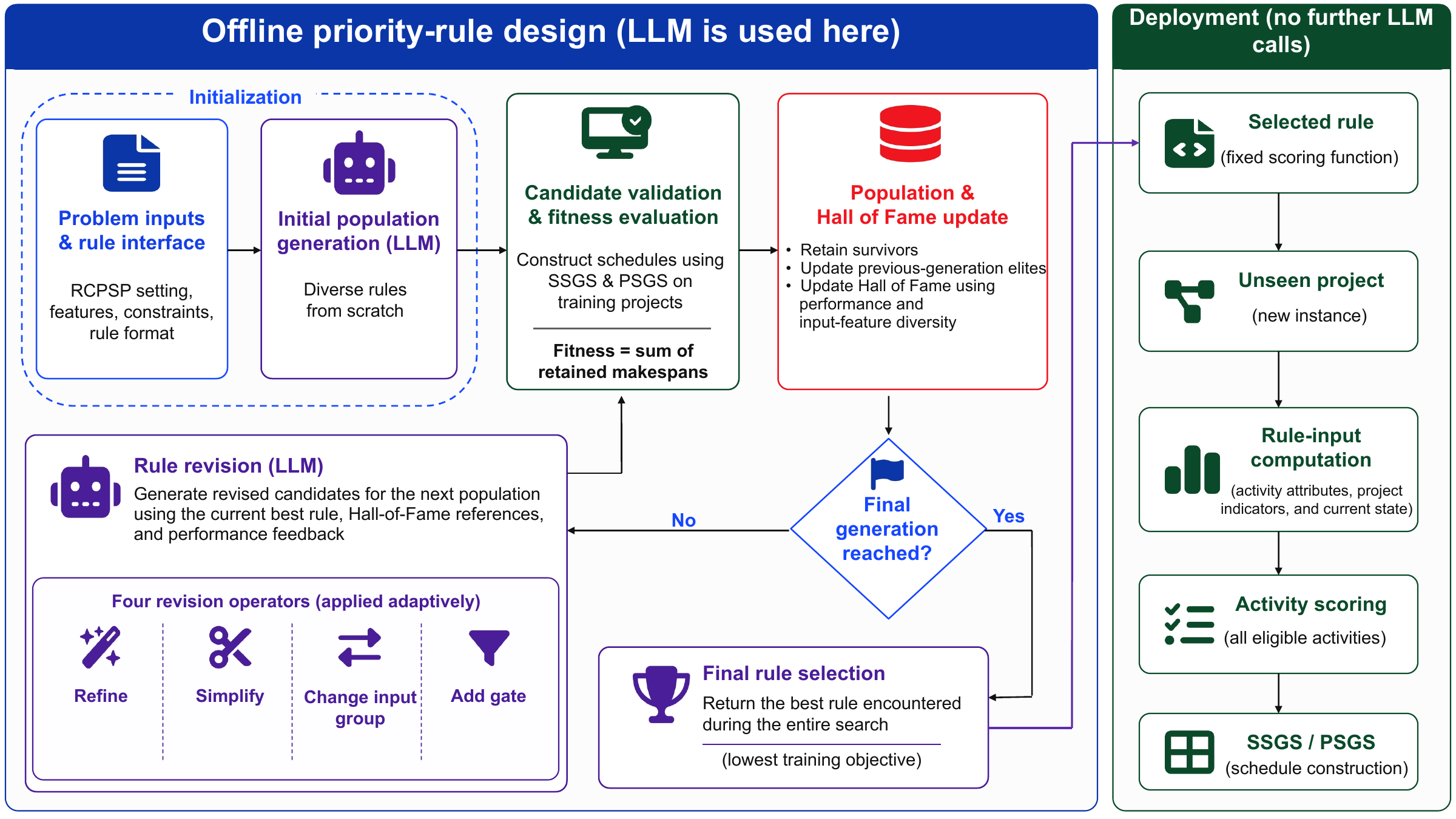}
\caption{Overview of the LLM-guided priority-rule design framework.}
\label{fig:framework_overview}
\end{figure}

The proposed framework embeds an LLM within a population-based search that follows an evolutionary update process. During offline rule design, the LLM generates and revises candidate rules through structured prompts, while candidate validation, fitness evaluation, and population management are performed outside the LLM. Candidate fitness is determined from schedules constructed on the training projects using both SSGS and PSGS. After the search, the rule with the lowest training objective encountered is returned and applied to new projects without further LLM calls. \Cref{fig:framework_overview} summarizes the overall workflow.

\subsection{Framework overview}
\label{sec:framework_overview}

As shown in the initialization part of \Cref{fig:framework_overview}, each search run begins by providing the LLM with the problem inputs and the rule interface. The problem inputs consist of the RCPSP setting and the available features. The RCPSP setting provides the problem background needed for rule design and explains how a priority rule is used during schedule construction. The available features serve as the rule inputs defined in \Cref{sec:priority_rules_sgs} and comprise activity-level attributes, project-level indicators, and schedule-construction information. The rule interface specifies the constraints imposed on the generated rule and its required output format. In particular, each rule must return a finite priority score for every eligible activity using only the permitted rule inputs and expression structures, and the activity with the lowest score is selected. The prompt also encourages compact and interpretable rules. Using this information, the LLM generates the initial population of priority rules from scratch.

Every newly generated candidate, during both initialization and subsequent generations, is produced through a separate LLM call. The framework does not modify existing formulas through hard-coded algebraic transformations. Instead, the initialization templates and revision operations are expressed as prompt instructions, and each call returns one candidate priority rule. The generated rules then enter candidate validation, where each output is parsed, checked against the prescribed interface and code restrictions, and deduplicated before fitness evaluation.

Each valid rule is evaluated by constructing schedules for every training project using both SSGS and PSGS. For each project, the lower of the two makespans is retained, and the sum of the retained makespans across the training set defines the rule's fitness. These fitness values initialize the search records by identifying the best rule encountered so far, selecting the generation elites, and constructing a Hall of Fame that balances training performance and diversity in the rule inputs used.

After the initial population has been evaluated, the framework updates the population records used to guide the search. Candidate rules are ranked by fitness, the best rule encountered so far is recorded, and the leading candidates are retained as generation elites. The Hall of Fame is constructed by considering both training performance and diversity in the sets of rule inputs used by the candidates. Selected rules from the Hall of Fame and the generation elites are retained as survivors for the next generation.

In subsequent generations, the best-performing rule in the Hall of Fame serves as the base rule for revision. The LLM receives this rule together with selected Hall-of-Fame references, recent performance feedback, and one assigned revision operation: refine, simplify, change input group, or add gate. The four revision operations are described in \Cref{sec:rule_revision}. It then revises the base rule according to this instruction and returns a new candidate priority rule. Each resulting candidate undergoes the same validation and fitness-evaluation procedure as the initial population, including schedule construction with both SSGS and PSGS. The revised candidates are combined with the retained survivors to form the next population, after which the generation elites, search records, and Hall of Fame are updated. This cycle repeats until the final generation is reached. The complete prompts are provided in the supplementary material.

At the end of the run, the rule with the lowest training objective encountered during the search is returned as the final rule. For a new project, the required rule inputs are computed and supplied to the selected priority function, which ranks eligible activities during schedule construction. As with traditional priority rules, the selected rule is used directly within SSGS and PSGS, and the schedule with the lower makespan is retained. Deployment requires neither a new rule-design search nor further LLM calls.

\subsection{Priority-rule representation and initial population}
\label{sec:rule_representation}

Let \(t\) index the decision points encountered during schedule construction, \(r\) denote a candidate priority rule, and \(f_r\) denote its executable scoring function. At decision point \(t\), rule \(r\) assigns eligible activity \(i\in E_t\) the score
\begin{equation}
q_{it}=f_r(\mathbf{x}_i,\mathbf{z}_t),
\label{eq:priority_function}
\end{equation}
where \(\mathbf{x}_i\) contains the activity-level inputs of activity \(i\), while \(\mathbf{z}_t\) contains the project-level and schedule-construction inputs shared by all activities eligible at decision point \(t\). Rule \(r\) remains fixed throughout schedule construction, but \(q_{it}\) may change across decision points because the schedule-construction inputs in \(\mathbf{z}_t\) are updated as the schedule is built.

At decision point \(t\), rule \(r\) selects
\begin{equation}
i_t^*=\min\operatorname*{arg\,min}_{i\in E_t}q_{it},
\label{eq:priority_rule_selection}
\end{equation}
so that the activity with the lowest score is selected. If several activities receive the same score, the activity with the smallest index is chosen.

The activity-level rule inputs are summarized in \Cref{tab:activity_inputs}. The symbols \(d_i\), \(r_{i,k}\), \(a_k\), \(S_i\), \(IP_i\), and \(N\) are defined in \Cref{sec:problem_description}. Let \(TS_i\) and \(TP_i\) denote the sets of total successors and total predecessors of activity \(i\), respectively, and let \(\ell_i\) denote the length of the longest path from activity \(i\) to the dummy finish activity \(N+1\). The resource equivalent duration \(RED_i\) combines the duration of activity \(i\) with its relative demand for renewable resources, following the original definition in \citet{cooper1976heuristics}.

The second column of \Cref{tab:activity_inputs} gives the field name used in the executable priority function, while the final column reports the corresponding definition before scaling. In later mathematical expressions, lower-case roman notation denotes the scaled field supplied to the rule; for example, \(\mathrm{lf}_i\) denotes the scaled value of \texttt{lf}. All activity-level rule inputs are nonnegative and scaled to \([0,1]\) before being supplied to the priority function.

\begingroup
\scriptsize
\setlength{\tabcolsep}{2.5pt}
\renewcommand{\arraystretch}{1.08}
\begin{longtable}{@{}>{\raggedright\arraybackslash}p{1.30cm}>{\raggedright\arraybackslash}p{2.35cm}>{\raggedright\arraybackslash}p{5.60cm}>{\raggedright\arraybackslash}p{7.05cm}@{}}
\caption{Activity-level inputs available to the priority rules.}
\label{tab:activity_inputs}\\
\toprule
Category & Rule input & Attribute & Definition before scaling \\
\midrule
\endfirsthead
\multicolumn{4}{l}{\textit{\Cref{tab:activity_inputs} continued}}\\
\toprule
Category & Rule input & Attribute & Definition before scaling \\
\midrule
\endhead
\bottomrule
\endfoot

Activity & \texttt{pt} & Activity processing time  & \(d_i\) \\

\addlinespace
Scheduling & \texttt{es} & Earliest start time  & \(ES_i\) \\
& \texttt{ef} & Earliest finish time  & \(EF_i\) \\
& \texttt{ls} & Latest start time  & \(LS_i\) \\
& \texttt{lf} & Latest finish time  & \(LF_i\) \\
& \texttt{slk} & Activity slack  & \(LS_i-ES_i\) \\

\addlinespace
Network & \texttt{rpw} & Ranked position weight  & \(d_i+\sum_{j\in IP_i}d_j\) \\
& \texttt{mts} & Number of total successors & \(\lvert TS_i\rvert\) \\
& \texttt{mtp} & Number of total predecessors  & \(\lvert TP_i\rvert\) \\
& \texttt{lfs} & Float of successors  & \((LS_i-ES_i)/\lvert TS_i\rvert\) \\
& \texttt{nrj} & Number of non-related activities  & \(N-\lvert TP_i\rvert-\lvert TS_i\rvert-1\) \\
& \texttt{lpf} & Length of the longest path following \(i\)  &  \(\ell_i\) \\
& \texttt{mis} & Number of immediate successors & \(\lvert S_i\rvert\) \\
& \texttt{msl} & Number of successors per level  &  \(\lvert TS_i\rvert/\ell_i\) \\
& \texttt{mtspt} & Total successor processing time  & \(d_i+\sum_{j\in TS_i}d_j\) \\

\addlinespace
Resource & \texttt{red} & Resource equivalent duration  & \(RED_i\) \\
& \texttt{cumred} & Cumulative resource equivalent duration & \(RED_i+\sum_{j\in S_i}RED_j\) \\
& \texttt{crwc} & Cumulative resource work content & \(\sum_{j\in S_i}\left(d_j\sum_{k=1}^{\lvert R\rvert}r_{j,k}\right)\) \\
& \texttt{trs} & Resource scarcity & \(\sum_{k=1}^{\lvert R\rvert}\dfrac{r_{i,k}}{a_k}\) \\
& \texttt{gres} & Resource requirement  & \(\sum_{k=1}^{\lvert R\rvert}r_{i,k}\) \\
& \texttt{grd} & Resource demand & \(d_i\sum_{k=1}^{\lvert R\rvert}r_{i,k}\) \\
& \texttt{wacru} & Weighted activity criticality  & \(w\sum_{j\in S_i}\dfrac{1}{1+d_j}+(1-w)\sum_{k=1}^{\lvert R\rvert}\dfrac{r_{i,k}}{a_k}\), with \(w=0.5\) \\
& \texttt{wrup} & Weighted resource utilization and precedence  & \(w\lvert S_i\rvert+(1-w)\sum_{k=1}^{\lvert R\rvert}\dfrac{r_{i,k}}{a_k}\) \\

\addlinespace
Other & \texttt{pop} & Product of selected priorities & \(d_i
   \times \lvert TS_i\rvert
   \times \left(\dfrac{\lvert TS_i\rvert}{\ell_i}\right)
   \times \left(d_i+\sum_{j\in IP_i}d_j\right)\) \\

\end{longtable}
\endgroup

In addition to the activity-level inputs in \Cref{tab:activity_inputs}, the priority rules may use the project-level indicators and schedule-construction inputs summarized in \Cref{tab:state_inputs}. The project-level indicators are computed once for each project, whereas the schedule-construction inputs are updated as SSGS or PSGS proceeds.

\begin{table}[htbp]
\centering
\caption{Project-level and schedule-construction rule inputs.}
\label{tab:state_inputs}

\begingroup
\scriptsize
\setlength{\tabcolsep}{3pt}
\renewcommand{\arraystretch}{1.12}

\begin{tabularx}{\textwidth}{
@{}
>{\raggedright\arraybackslash}p{3.05cm}
@{\hspace{0.40cm}}
>{\raggedright\arraybackslash}p{3.15cm}
@{\hspace{0.15cm}}
>{\raggedright\arraybackslash}X
>{\centering\arraybackslash}p{1.15cm}
@{}}
\toprule
Category & Rule input & Interpretation & Range \\
\midrule

\makecell[l]{Schedule-construction\\state}
& \texttt{progress}
& Fraction of activities scheduled in SSGS or completed in PSGS
& \([0,1]\) \\

& \texttt{queue\_length}
& Number of currently eligible activities, \(\lvert E_t\rvert\)
& \(\{0,\ldots,N\}\) \\

& \makecell[l]{\texttt{avg\_res\_}\\\texttt{utilization}}
& Resource utilization of the partial schedule: cumulative in SSGS and instantaneous in PSGS
& \([0,1]\) \\

\addlinespace
\makecell[l]{Project-network\\indicators}
& \texttt{sp}
& Serial/parallel indicator; lower values represent more parallel networks and higher values more serial networks
& \([0,1]\) \\

& \texttt{ad}
& Activity distribution; variation in the distribution of activities across network levels
& \([0,1]\) \\

& \texttt{la}
& Length of arcs; prevalence of arcs connecting adjacent network levels
& \([0,1]\) \\

& \texttt{tf}
& Topological float ; flexibility in assigning activities to network levels
& \([0,1]\) \\

\addlinespace
\makecell[l]{Project-resource\\indicators}
& \texttt{rc}
& Resource constrainedness; higher values indicate tighter resource conditions
& \([0,1]\) \\

& \texttt{rs}
& Resource strength; higher values indicate greater available capacity relative to demand
& \([0,1]\) \\

& \texttt{rf}
& Resource factor; proportion of activity--resource pairs with positive requirements
& \([0,1]\) \\

& \texttt{ru}
& Resource use; average number of resource types required per activity
& \([0,\lvert R\rvert]\) \\

\bottomrule
\end{tabularx}

\endgroup
\end{table}

The \texttt{progress} input is defined separately for the two schedule generation schemes:
\begin{equation}
\mathrm{progress}^{\mathrm{SSGS}}_t
=
\frac{\lvert A^{\mathrm{scheduled}}_t\rvert}{N},
\qquad
\mathrm{progress}^{\mathrm{PSGS}}_t
=
\frac{\lvert A^{\mathrm{completed}}_t\rvert}{N},
\label{eq:construction_progress}
\end{equation}
where \(A^{\mathrm{scheduled}}_t\) is the set of non-dummy activities already inserted into the partial schedule before decision point \(t\), and \(A^{\mathrm{completed}}_t\) is the set of non-dummy activities completed by the scheduling time associated with decision point \(t\). In both schemes, \(\mathrm{progress}_t\) ranges from 0 to 1 and takes the same value for all eligible activities at a given decision point. It allows a priority rule to change its decision logic as schedule construction advances.

Together, \Cref{tab:activity_inputs,tab:state_inputs} define the rule inputs available to each candidate. The executable priority function may combine these inputs with numerical constants using arithmetic operations (addition, subtraction, multiplication, and division), comparisons, conditional branches, and a predefined set of elementary mathematical functions.

Before fitness evaluation, each LLM-generated function is checked for compliance with the required priority-function interface and code restrictions. A candidate is rejected if it uses unsupported rule inputs or functions, external code or data, randomness, or other prohibited program structures. The complete validation rules are provided in the supplementary materials.

The prompt descriptions do more than enumerate the available rule inputs. They incorporate RCPSP domain knowledge by explaining the scheduling roles of these inputs, indicating their expected priority directions, identifying selected inputs as strong priority-rule anchors, and warning against redundant input combinations. This information guides rule generation without prescribing a fixed formula.

The initial population is generated through independent LLM calls, with each call asking the LLM to design a priority rule from scratch. Each call receives the common problem inputs, rule interface, and code constraints described above. To encourage structural variation, the initialization prompt presents four alternative starting structures:

\begin{enumerate}
\item a rule without conditional branches, centered on one timing-related activity attribute, such as latest start, latest finish, or slack, together with one additional term representing successor-related information or resource requirements;
\item a rule with one condition based on schedule-construction progress and different priority expressions before and after a progress threshold;
\item a rule with one condition based on current average resource utilization and different priority expressions below and above a utilization threshold; and
\item a rule with exactly one conditional statement, with one branch emphasizing timing-related or successor-related information and the other emphasizing resource demand or scarcity.
\end{enumerate}

The four structures are alternative starting templates and are not combined in a single rule. For example, a progress-gated rule may use timing and successor-related terms before a progress threshold and resource-related terms after it, thereby illustrating both the second structure and the branch specialization described in the fourth. A utilization-gated rule follows a similar form but uses current average resource utilization as the switching condition. A complete LLM-generated rule with progress-dependent switching is presented in \Cref{eq:reference_rule}.

These structures guide the initial form of a rule, while the LLM determines the specific inputs, coefficients, thresholds, and mathematical expression. No previously generated rules, Hall-of-Fame references, or performance feedback are supplied during initialization. After validation and deduplication, the accepted rules form the initial population. These initialization structures differ from the revision operations applied to existing rules in subsequent generations, which are described in \Cref{sec:rule_revision}.

\subsection{Fitness evaluation and performance feedback}
\label{sec:fitness_feedback}

For a candidate rule \(r\) and training set \(\mathcal{T}\), one schedule is constructed with SSGS and another with PSGS for each project \(p\in\mathcal{T}\). The lower of the two makespans is retained, and the fitness of the rule is defined as
\begin{equation}
F(r)
=
\sum_{p\in\mathcal{T}}
\min\left\{
C_{rp}^{\mathrm{SSGS}},
C_{rp}^{\mathrm{PSGS}}
\right\},
\label{eq:training_fitness}
\end{equation}
where \(C_{rp}^{\mathrm{SSGS}}\) and \(C_{rp}^{\mathrm{PSGS}}\) are the makespans obtained by applying rule \(r\) to project \(p\) under the two schedule generation schemes. Because \(F(r)\) is the sum of the retained makespans, a lower fitness value corresponds to better schedule quality. Since all candidates in a search run are evaluated on the same training set, their fitness values can be used directly to rank the candidate rules.

After each generation is evaluated, the results are converted into summary feedback for the next revision prompts. The LLM does not receive the individual schedules or instance-level makespans. Instead, each prompt reports the best fitness in the previous generation, the best fitness found during the search, the average population fitness, the number of consecutive generations without a meaningful improvement, and a numerical improvement target. Each Hall-of-Fame reference rule is also accompanied by its fitness. These values are calculated after the candidate rules are evaluated on the training projects; the LLM only generates and revises candidate rules. \Cref{sec:rule_revision} explains how the revision prompt combines this feedback with the base rule, reference rules, and assigned revision operation.

\subsection{Rule revision, population update, and final rule selection}
\label{sec:rule_revision}

The initialization and revision prompts serve different purposes. The initialization prompt asks the LLM to design a candidate rule from scratch, whereas the revision prompt asks it to modify an existing rule. \Cref{tab:prompt_structure} compares the information supplied in the two cases.

From the second generation onward, the best-performing rule in the Hall of Fame serves as the base rule for the new candidates in that generation. Each candidate is generated through a separate LLM call. In addition to the base rule, the revision prompt provides prior priority-rule knowledge, selected Hall-of-Fame references, recent performance feedback, one assigned revision operation, and the required output format. The LLM then revises the base rule according to this information and returns one candidate priority rule.

\begin{table}[htbp]
\centering
\caption{Information supplied to the LLM during initialization and rule revision.}
\label{tab:prompt_structure}
\begingroup
\footnotesize
\setlength{\tabcolsep}{4pt}
\renewcommand{\arraystretch}{1.10}
\begin{tabularx}{\textwidth}{@{}p{3.2cm}>{\raggedright\arraybackslash}X>{\raggedright\arraybackslash}X@{}}
\toprule
Prompt element & Initial population & Subsequent generations \\
\midrule
Starting point & No base rule; each candidate is designed from scratch & Current best Hall-of-Fame rule, supplied as the base rule to revise \\
Common rule-design specification & RCPSP setting, available rule inputs, minimum-score convention, rule interface, and code constraints & Same common specification \\
Prior priority-rule knowledge
& Guidance on input meanings, expected priority directions, strong anchors, and redundant input combinations
& Same guidance \\
Within-run search information & None, because no previous population has been evaluated & Selected Hall-of-Fame reference rules and their fitness values, recent aggregate performance feedback, and summaries of the rule inputs used by the reference rules \\
Generation instruction & Prompt presenting four alternative starting structures & One assigned revision operation: refine, simplify, change input group, or add gate \\
Required output & One candidate priority rule & One revised candidate priority rule \\
\bottomrule
\end{tabularx}
\endgroup
\end{table}

The requested revision is guided by one of four prompt-based operations:

\begin{enumerate}[label=(\roman*),leftmargin=2.4em,itemsep=0.25em,topsep=0.35em,parsep=0pt]
\item \emph{Refine}: adjusts one or two important parts of the base rule by strengthening a useful term, changing a coefficient, or removing a weak term, without rescaling the entire expression.
\item \emph{Simplify}: removes redundant structure and seeks a more compact rule centered on one principal ranking criterion, with at most two additional terms that adjust this criterion.
\item \emph{Change input group}: changes the group of inputs on which the rule mainly relies, for example by shifting among the scheduling, network, and resource attributes in \Cref{tab:activity_inputs}, or toward current average resource utilization in \Cref{tab:state_inputs}.
\item \emph{Add gate}: introduces one conditional statement that switches between two priority-score formulas. The gate is based on either schedule-construction progress or current average resource utilization. Fixed project-network and project-resource indicators may still appear within either formula.
\end{enumerate}

For illustration, consider the base rule
\(q_{it}=\mathrm{lf}_i-0.5\,\mathrm{mtspt}_i-0.2\,\mathrm{grd}_i\).
A refine prompt may change the coefficient of \(\mathrm{grd}_i\), while a simplify prompt may remove this term. A change-input-group prompt may replace the timing-based main criterion with a resource-based criterion. An add-gate prompt may use the base expression before a progress threshold and a different expression afterward. These examples illustrate the requested direction of revision; the LLM determines the actual revised rule.

Before describing the population update, we distinguish the roles of the generation elites, Hall of Fame, base and reference rules, survivors, and the best rule found over the complete run. Generation elites are the best-performing rules in the current generation. The Hall of Fame is a persistent archive updated across generations using both fitness and diversity in rule-input use. Its best-performing member serves as the base rule, while selected Hall-of-Fame members are supplied as reference rules in the revision prompts. Survivors are selected Hall-of-Fame members and generation elites carried directly into the next population. The best rule found over the complete run is tracked separately and returned at the end.

After all rules in a generation have been evaluated, they are ranked by fitness, and approximately the top 15\% of the population are designated as generation elites. The Hall of Fame has a capacity of approximately 25\% of the population size. To update it, the existing Hall-of-Fame members are combined with the rules in the current generation, and duplicates are removed. One quarter of the archive capacity is filled with the rules having the lowest fitness values. The remaining positions are selected sequentially using a criterion that favors both lower fitness and greater diversity in rule-input use. For two rule-input sets \(A\) and \(B\), the Jaccard distance is
\(1-\lvert A\cap B\rvert/\lvert A\cup B\rvert\)
\citep{levandowsky1971distance}. For each candidate, diversity is measured by the minimum of this distance over the rules already selected into the Hall of Fame. After the update, the Hall of Fame is sorted by fitness, and its best-performing members are supplied as reference rules in subsequent revision prompts. Because Hall-of-Fame membership is diversity-aware, the resulting reference set is drawn from an archive that reflects both fitness and rule-input diversity.

Approximately 15\% of the population is carried into the next generation as distinct survivors. The survivor quota is divided equally between rules selected from the updated Hall of Fame and elites from the most recent generation. If a rule is selected from both sources, it is retained only once. Survivor positions left unfilled because of overlap or an insufficient number of distinct rules are reassigned to newly generated candidates. The exact archive, elite, and survivor sizes used in the experiments are reported in \Cref{sec:evaluation_protocol}.

After survivors have been selected, the number of new candidates required to complete the next population is determined. The corresponding LLM calls are then allocated across the four revision operations according to recent search progress, as shown in \Cref{tab:revision_allocation}. Each call receives the common base rule and one assigned revision operation and returns one new candidate rule. The approximate proportions are converted into integer call counts, and the operation assignments are randomly shuffled before the calls are submitted.

\begin{table}[htbp]
\centering
\caption{Approximate allocation of LLM calls across revision operations.}
\label{tab:revision_allocation}
\begingroup
\footnotesize
\setlength{\tabcolsep}{5pt}
\renewcommand{\arraystretch}{1.10}
\begin{tabularx}{\textwidth}{@{}l*{4}{>{\centering\arraybackslash}X}@{}}
\toprule
Search status & Refine & Simplify & \makecell{Change input\\group} & Add gate \\
\midrule
Normal progress & 50\% & 20\% & 25\% & 5\% \\
Stagnation & 30\% & 20\% & 30\% & 20\% \\
\bottomrule
\end{tabularx}
\endgroup
\end{table}

The search is classified as stagnant when the best Hall-of-Fame fitness values vary by less than a predefined threshold over a fixed number of recent generations. Under stagnation, the proportion assigned to refinement is reduced, while larger proportions are assigned to changing the main input group and adding conditional structures.

Sampling temperature provides an additional means of controlling variation in the generated rules. An LLM generates a rule one token at a time by assigning probabilities to possible next tokens and sampling from this distribution. The temperature controls how concentrated the distribution is. At a temperature of 1, the original probability distribution is used. Lower temperatures place more weight on high-probability tokens and reduce variation, whereas higher temperatures make lower-probability alternatives more likely. 

The temperatures are kept close to 1 so that the assigned revision operation determines the main direction of the requested change, while temperature makes only a moderate adjustment to output diversity. Refine and simplify use lower temperatures than change input group and add gate, and the temperatures are increased moderately under stagnation. The numerical temperature values, stagnation window, and threshold used in the experiments are reported in \Cref{sec:evaluation_protocol}.

For clarity, we use the term \textit{search components} for the main elements of the framework described above. \textit{Cross-generation dependence} means that later generations use rules or information obtained in earlier generations. \textit{Base-rule passing} supplies a previously identified rule for revision, while \textit{performance feedback} supplies aggregate fitness information. \textit{Operator-specific prompts} refer to the four revision operations. \textit{Hall-of-Fame references} are selected archive members supplied during revision, while \textit{diversity-aware Hall-of-Fame management} uses rule-input diversity when updating the archive from which those references are drawn. \textit{Survivor retention} carries previously evaluated rules directly into the next population.

The survivors and newly generated candidates together form the next population. After the new candidates have been evaluated, the generation elites, Hall of Fame, performance records, and run-level best rule are updated, and the same cycle begins for the following generation. The process continues until the final generation is reached. The rule with the lowest training fitness encountered over the complete run is then returned as the final rule, regardless of the generation in which it was generated. The complete initialization and revision prompts are provided in the supplementary materials.

\section{Experimental design}
\label{sec:experimental_design}

This section introduces the benchmark data and experimental settings. \Cref{sec:benchmark_design} describes the datasets and their roles in training, testing, and post-selection evaluation, while \Cref{sec:evaluation_protocol} presents the evaluation protocol and implementation details.

\subsection{Benchmark data and training--test design}
\label{sec:benchmark_design}

The experiments use instances from the PSPLIB and RanGen benchmark families, together with three datasets containing projects with more than 1,000 activities. \Cref{tab:dataset_roles} summarizes the size and experimental role of each dataset.

\begin{table}[htbp]
\centering
\caption{Benchmark datasets and experimental roles.}
\label{tab:dataset_roles}
\scriptsize
\begin{tabularx}{\textwidth}{@{}llrr>{\raggedright\arraybackslash}X@{}}
\toprule
Dataset group & Dataset & Instances & Activities & Experimental role \\
\midrule
PSPLIB & J30 & 480 & 30 & Training data \\
& J60 & 480 & 60 & Test data \\
& J90 & 480 & 90 & Test data \\
& J120 & 600 & 120 & Test data \\
\addlinespace
RanGen & RG30 & 1,800 & 30 & Test data for J30-trained rules; training data for J30+RG30-trained rules \\
& RG300 & 480 & 300 & Test data \\
\addlinespace
Large projects & Multi-projects & 226 & 1,440 & Post-selection evaluation \\
& LPP & 80 & 1,800 & Post-selection evaluation \\
& LPSP & 81 & 1,800 & Post-selection evaluation \\
\bottomrule
\end{tabularx}
\end{table}

PSPLIB comprises the J30, J60, J90, and J120 sets, whose instances contain 30, 60, 90, and 120 activities, respectively \citep{kolisch1997psplib}. J30, J60, and J90 each contain 480 instances generated from 48 combinations of three ProGen input parameters: network complexity (NC), resource factor (RF), and resource strength (RS), with ten instances for each combination. NC controls the average number of non-redundant precedence arcs per activity, while RF and RS correspond to the resource factor and resource strength indicators introduced in \Cref{tab:state_inputs}. J120 contains 600 instances generated from 60 combinations of the same parameters. The four sets provide projects of increasing size within the common PSPLIB generation framework.

RG30 was generated using RanGen and contains 1,800 instances divided into five subsets \citep{demeulemeester2003rangen,vanhoucke2008evaluation}. Compared with J30, RG30 covers a substantially broader range of network structures. For example, the serial/parallel indicator (SP) ranges from 0.10 to 0.90 in RG30, compared with 0.17 to 0.41 in J30, while the coefficient of network complexity (CNC) ranges from 0.23 to 6.47, compared with 1.40 to 2.07. Similar differences are observed for other network-topology indicators. RG300 contains 480 RanGen-generated instances with 300 activities and is used to evaluate whether rules trained on 30-activity instances generalize to projects with 300 activities \citep{debels2007decomposition}.

The experiments compare two training datasets for priority-rule design: J30 alone and the expanded J30+RG30 set. When rules are trained on J30, all 480 J30 instances are used for training, while J60, J90, J120, RG30, and RG300 are reserved for testing. When rules are trained on J30+RG30, all 1,800 RG30 instances are added to J30, yielding 2,280 training instances, while J60, J90, J120, and RG300 remain as test sets. Thus, RG30 is unseen during training in the first case but contributes directly to training in the second. Given its broader coverage of network-topology indicators, adding RG30 also broadens the range of project structures represented in the training data. All other components of the LLM-guided search remain unchanged.

Three additional datasets are used to evaluate performance on projects with more than 1,000 activities. The Multi-projects dataset was introduced by \citet{VanEynde2020a} and originally contains 2,254 instances with 1,440 activities each. Following \citet{luo2022efficient}, every tenth instance is retained, resulting in 226 evaluation instances. These instances have highly parallel network structures, with a mean SP value of 0.03 and a range from 0.01 to 0.04.

The two large single-project datasets, LPP and LPSP, were constructed by \citet{luo2022efficient}. For LPP, a RanGen network containing six non-dummy activities defines the upper-level structure, and each activity is replaced by one RG300 project. Grouping the 480 RG300 instances into sets of six produces 80 large projects with \(6\times300=1{,}800\) activities each. The resulting projects are highly parallel, with a mean SP value of 0.018 and a range from 0.01 to 0.03. Each LPSP instance is formed by combining 36 component projects with 50 activities each, again resulting in 1,800 activities. SP takes nine values from 0.10 to 0.90. At each SP level, one instance is generated for each of nine resource-constrainedness levels from 0.10 to 0.90, yielding 81 instances. LPSP therefore covers a broader range of network structures than Multi-projects and LPP.

The three large-project datasets are used exclusively for the post-selection evaluation reported in \Cref{sec:results_large_projects}; none of their instances or results contribute to rule design or rule selection.

\subsection{Evaluation protocol and experimental settings}
\label{sec:evaluation_protocol}

This subsection presents the evaluation protocol and the main implementation and statistical settings.

The rules returned by the search runs are evaluated on the test sets using the same best-of-SSGS/PSGS procedure used to compute candidate fitness in \Cref{sec:fitness_feedback}. For each rule and instance, schedules are constructed with both schemes, and the lower makespan is retained. This protocol follows \citet{luo2026automated}, which shows that combining the two schemes improves solution quality with only marginal additional computational effort. The schedule-generation and tie-breaking procedures are deterministic, so repeated evaluation of the same rule on the same instance produces the same retained makespan.

Performance is calculated separately for each test set using the average percentage deviation from the best-known lower bound:
\begin{equation}
\operatorname{AvgDevLB}(r,\mathcal{I})
=
\frac{100}{\lvert\mathcal{I}\rvert}
\sum_{p\in\mathcal{I}}
\frac{C_{rp}-LB_p}{LB_p},
\label{eq:avgdevlb}
\end{equation}
where \(\mathcal I\) denotes the test set under consideration in \Cref{tab:dataset_roles}, \(p\in\mathcal I\) indexes a project instance, and
\(C_{rp}=\min\{C_{rp}^{\mathrm{SSGS}},C_{rp}^{\mathrm{PSGS}}\}\)
is the retained makespan obtained by rule \(r\) on that instance.
\(LB_p\) is the best-known lower bound for project instance \(p\), as reported by \citet{vanhoucke2018tool}. Lower AvgDevLB values indicate better performance.

Because lower-bound data are unavailable for the three large-project datasets, performance on these datasets is measured by the average percentage deviation from the best-known solution:
\begin{equation}
\operatorname{AvgDevBKS}(r,\mathcal{I})
=
\frac{100}{\lvert\mathcal{I}\rvert}
\sum_{p\in\mathcal{I}}
\frac{C_{rp}-BKS_p}{BKS_p},
\label{eq:avgdevbks}
\end{equation}
where \(BKS_p\) is the best-known makespan for project instance \(p\), as reported by \citet{luo2022efficient}. Lower AvgDevBKS values indicate better performance.

Before conducting the main experiments, we screened four combinations of population size and generation count. Each combination was evaluated over three independent runs using only the training objective on J30. Test-set performance was not used for configuration selection, and all other framework settings were identical across the four configurations. \Cref{tab:configuration_screening} reports the results.

\begin{table}[htbp]
\centering
\caption{Configuration screening on the J30 training set.}
\label{tab:configuration_screening}
\scriptsize
\begin{tabular}{@{}llllll@{}}
\toprule
\makecell{Population\\size} & Generations & Run 1 & Run 2 & Run 3 & Mean \\
\midrule
20 & 15 & 29,271 & 29,242 & 29,253 & 29,255 \\
20 & 25 & 29,238 & 29,233 & 29,214 & 29,228 \\
50 & 25 & 29,200 & 29,215 & 29,213 & \textbf{29,209} \\
100 & 15 & 29,217 & 29,228 & 29,186 & 29,210 \\
\bottomrule
\end{tabular}
\end{table}

As shown in \Cref{tab:configuration_screening}, the configuration with a population size of 50 and 25 generations achieved the lowest mean training objective. Its range across the three runs was also smaller than that of the nearly tied configuration with a population size of 100 and 15 generations. It was therefore used in all subsequent LLM-guided experiments.

For the selected population size of 50, the proportional population-update settings in \Cref{sec:rule_revision} correspond to the following values. The generation-elite set contains up to eight rules, approximately 15\% of the population, and the Hall of Fame contains up to 12 rules, approximately 25\% of the population. Three positions, corresponding to 25\% of the Hall-of-Fame capacity, are filled solely according to fitness, while the remaining positions are selected using the fitness-and-diversity criterion. Up to eight rules, approximately 15\% of the population, are retained as survivors, with at most four selected from the Hall of Fame and four from the preceding generation's elites.

Following common practice in LLM-based automated heuristic design, we compare search effort using a nominal new-candidate evaluation budget, which counts newly generated candidates and excludes retained survivors \citep{zhang2024understanding}. When eight distinct survivors are retained, the initial generation contains 50 new candidates and each of the remaining 24 generations contains 42, giving a nominal new-candidate evaluation budget of \(50+24\times42=1{,}058\) per search run. This budget measures the number of newly evaluated candidate rules rather than the total number of schedules constructed during fitness evaluation. The former reflects how many rules are explored by the search, whereas the latter also depends on the training-set size and the procedure used to evaluate each rule \citep{liu2024evolution,zhang2024understanding}. We therefore use the number of newly evaluated candidate rules as the common measure of search effort and examine computational time separately in \Cref{sec:results_baselines}.

The initial population is generated with a sampling temperature of 1.00. From the second generation onward, up to eight survivors are retained from the previous population. When all eight survivor positions are filled, 42 of the 50 population positions are therefore filled by newly generated candidates. \Cref{tab:revision_settings} reports the allocation of these 42 LLM calls across the four revision operations, together with the sampling temperature used for each operation. Each entry gives the number of LLM calls, followed by the sampling temperature in parentheses.

\begin{table}[htbp]
\centering
\caption{Revision-operation settings for a standard generation of 42 newly generated candidates.}
\label{tab:revision_settings}
\begingroup
\footnotesize
\setlength{\tabcolsep}{5pt}
\renewcommand{\arraystretch}{1.10}
\begin{tabularx}{\textwidth}{@{}l*{4}{>{\centering\arraybackslash}X}@{}}
\toprule
Search status & Refine & Simplify & \makecell{Change input\\group} & Add gate \\
\midrule
Normal progress & 22 (0.95) & 8 (0.95) & 10 (1.08) & 2 (1.08) \\
Stagnation & 14 (1.02) & 8 (1.02) & 12 (1.12) & 8 (1.12) \\
\bottomrule
\end{tabularx}
\endgroup
\end{table}

The search is classified as stagnant when the range of the best Hall-of-Fame fitness values over the four most recent generations is below 12 units of the training objective.

The main LLM-guided configurations are each evaluated over ten independent search runs, and the reported results are means of the ten run-level performance values. The LLM backbone refers to the model used to generate and revise candidate rules. \Cref{sec:results_backbones} examines whether the performance of the framework depends on this model choice.

The search framework and the SSGS and PSGS procedures are implemented in Python 3.13. All experiments are conducted on a Windows 11 computer equipped with an Intel Core i9-12900H processor and 32 GB of RAM. Each generated rule is implemented as a function named \texttt{priority\_score}, with \texttt{activity} and \texttt{state} as its two arguments; the rule-validation procedure is described in \Cref{sec:rule_representation}. Gemini Flash\footnote{API model identifier: \texttt{gemini-3-flash-preview}.} is used as the default LLM backbone. An LLM call denotes one API request used to generate or revise a candidate rule. LLM-call budgets are reported separately where relevant.

Statistical comparisons are conducted separately for each test set. When both sides of a comparison are represented by multiple independent search runs, each run contributes one observation: the AvgDevLB obtained by the rule returned from that run on the complete test set. Thus, when both configurations are evaluated over ten runs, the comparison is based on two samples of ten run-level AvgDevLB values and is conducted using a Mann--Whitney U test. A one-sided test is used when a directional alternative hypothesis is specified in advance, whereas a two-sided test is used when no direction is prespecified.

When an LLM-guided configuration is compared with a deterministic single priority rule, we follow the instance-level comparison procedure used in \citet{luo2023automated,luo2026automated}. For each test instance, the percentage deviations obtained by the ten independently generated LLM rules are first averaged, and the resulting mean is paired with the deviation obtained by the deterministic rule on the same instance. A one-sided Wilcoxon signed-rank test is then applied to these instance-level pairs.

The direction of each alternative hypothesis and the comparison baseline are stated in the corresponding results subsection. All tests use a significance level of \(\alpha=0.05\). Across the result tables, a superscript asterisk denotes a statistically significant difference from the comparison baseline identified in the corresponding subsection.

\section{Computational results}
\label{sec:results}

This section evaluates the proposed framework from five perspectives: comparison with established priority-rule methods, the contribution of the search design, the effect of expanded training, robustness across LLM backbones, and performance on large projects. \Cref{tab:experiment_overview} summarizes the purpose and main result of each experiment before the detailed results are presented.

\begin{table}[htbp]
\centering
\caption{Overview of the computational experiments.}
\label{tab:experiment_overview}
\footnotesize
\setlength{\tabcolsep}{4pt}
\renewcommand{\arraystretch}{1.14}
\begin{tabularx}{\textwidth}{@{}>{\centering\arraybackslash}p{1.5cm}>{\raggedright\arraybackslash}X>{\raggedright\arraybackslash}X@{}}
\toprule
Subsection & Purpose & Main result \\
\midrule
\Cref{sec:results_baselines}
& Compare LLM-designed rules with traditional and GP-designed rules under different search efforts.
& The LLM-designed rules outperform the single-rule baselines and show favorable performance--search-effort trade-offs relative to GP. \\

\Cref{sec:results_ablation}
& Assess the contribution of the main search components.
& The full framework outperforms all five ablation variants. \\

\Cref{sec:results_expanded}
& Assess the effect of broader training data on rule performance.
& Broader training data improve performance across all test sets. \\

\Cref{sec:results_backbones}
& Assess robustness across LLM backbones and compare API usage and cost.
& Performance remains robust across backbones, although model choice affects rule quality and cost. \\

\Cref{sec:results_large_projects}
& Assess whether selected rules retain their performance on much larger projects.
& The selected rules remain effective on large projects, with particularly strong results on highly parallel datasets. \\
\bottomrule
\end{tabularx}
\end{table}

\subsection{Comparison with priority-rule baselines}
\label{sec:results_baselines}

We first compare the test-set performance of the priority rules generated by the LLM-guided framework with the performance of traditional and GP-designed rules. For the GP comparisons, we use the GPHH framework of \citet{luo2022efficient} under two configurations. The matched-budget configuration uses a population size of 50 and 25 generations, matching the LLM-guided framework. The original configuration of \citet{luo2022efficient} uses a population size of 1,000 and 50 generations.

The traditional baselines comprise worst-case slack (WCS, \citet{Kolisch1996}), latest finish time (LFT), SinglePR (the best single traditional rule for each test set), and AllPR (the best schedule obtained from all 36 traditional rules for each instance). WCS and LFT are individual priority rules. For each test set, SinglePR denotes the one rule among the 36 traditional priority rules examined in \citet{luo2022efficient,luo2026automated} that obtains the lowest AvgDevLB over the complete test set. The selected rule is applied to every instance in that set. AllPR instead applies all 36 rules to each instance and retains the best schedule obtained for that instance. AllPR is therefore treated as a multi-rule reference and excluded from single-rule comparisons.

All methods are evaluated under the common test protocol described in \Cref{sec:evaluation_protocol}. Both GP configurations retain the activity-attribute terminal set used by \citet{luo2022efficient}. The LLM rule interface includes the same activity attributes and additionally permits the project-level indicators and schedule-construction information described in \Cref{sec:rule_representation}. The LLM-guided search uses all 480 J30 instances as training data, whereas both GP configurations use the 96-instance J30 subset adopted by \citet{luo2022efficient}. This subset contains the first two instances from each of the 48 PSPLIB parameter combinations. \citet{luo2022efficient} compared complete and sampled training data drawn from both J30 and RG30 and showed that carefully sampled training sets could retain comparable test-set performance while substantially reducing the computational effort required for GP training. Because the LLM-guided and GP searches use training sets of different sizes, similar nominal new-candidate evaluation budgets do not imply the same number of schedule constructions during fitness evaluation. 

Using the nominal new-candidate evaluation budget defined in \Cref{sec:evaluation_protocol}, the LLM-guided framework has a nominal budget of 1,058 per run. The matched-budget GP evaluates approximately 1,030 newly generated candidate rules per run on average, with minor variation caused by its stopping rule. This is close to the LLM-guided framework's nominal budget of 1,058, so the two searches evaluate similar numbers of new candidate rules per run. In contrast, the original GP configuration evaluates approximately 50,000 newly generated candidate rules per run, representing the substantially larger search budget used by \citet{luo2022efficient}. The LLM-guided framework and both GP configurations are reported as means over ten independent runs, whereas the traditional priority-rule results are deterministic. Statistical comparisons follow the procedures defined in \Cref{sec:evaluation_protocol}. For comparisons with the deterministic single-rule baselines and the GP configurations, the alternative hypothesis is that the LLM-guided framework achieves lower AvgDevLB. Boldface identifies the lowest AvgDevLB among the single-rule approaches.

\begin{table}[htbp]
\centering
\caption{Test-set AvgDevLB (\%) against priority-rule baselines.}
\label{tab:baseline_comparison}
\scriptsize
\setlength{\tabcolsep}{3pt}
\renewcommand{\arraystretch}{1.08}
\sisetup{detect-weight=true,mode=text,table-number-alignment=center}
\begin{tabular*}{\textwidth}{@{}l@{\hspace{10pt}}l@{\extracolsep{\fill}}*{5}{S[table-format=2.3]}c@{}}
\toprule
Category & Rule/configuration & {J60} & {J90} & {J120} & {RG30} & {RG300} & \makecell{Nominal new-candidate\\evaluation budget} \\
\midrule
LLM-guided & This study & \bfseries 4.232 & \bfseries 4.599 & 11.743 & \bfseries 6.392 & \bfseries 12.203 & 1,058 \\
Traditional & WCS & 5.746\sig & 5.742\sig & 13.203\sig & 8.550\sig & 13.656\sig & N/A \\
& LFT & 4.555\sig & 4.947\sig & 12.221\sig & 8.098\sig & 13.170\sig & N/A \\
& SinglePR & 4.372\sig & 4.845\sig & 12.166\sig & 7.145\sig & 13.012\sig & N/A \\
& AllPR & 3.334 & 4.038 & 10.604 & 4.449 & 11.748 & N/A \\
GP-designed & Matched budget & 4.477\sig & 4.878\sig & 12.176\sig & 7.402\sig & 13.140\sig & 1,030 \\
& Original configuration & 4.250 & 4.620 & \bfseries 11.741 & 6.498\sig & 12.546\sig & $\approx$50,000 \\
\bottomrule
\end{tabular*}
\end{table}

\Cref{tab:baseline_comparison} shows that the LLM-guided framework outperforms WCS, LFT, and SinglePR on all five test sets, demonstrating a consistent advantage over individual traditional priority rules. AllPR achieves lower AvgDevLB than the LLM-guided framework on all five test sets. However, AllPR applies all 36 traditional priority rules to each instance and retains the best resulting schedule, whereas each run of the LLM-guided framework returns a single priority rule. AllPR therefore provides a demanding multi-rule reference. We revisit this comparison using the expanded training set and again on large projects in \Cref{sec:results_expanded,sec:results_large_projects}.

Under a comparable nominal new-candidate evaluation budget, the LLM-guided framework obtains lower mean AvgDevLB than the matched-budget GP on all five test sets, with all differences statistically significant at the 0.05 level. The original GP configuration, using an approximately 50-times larger budget, narrows the gap. The LLM-guided framework still obtains lower AvgDevLB on J60, J90, RG30, and RG300, while the results on J120 are nearly identical; the differences on RG30 and RG300 are statistically significant.

The preceding comparisons relate test-set performance to search effort measured by the nominal new-candidate evaluation budget. We next examine the wall-clock time of the LLM-guided search. Wall-clock time refers to the elapsed real time from the start to the end of a search run, rather than CPU time. It is used because the search combines remote LLM requests with parallel local candidate-rule evaluations. Local candidate-rule evaluations are performed in parallel, and LLM requests within a generation can also be submitted concurrently. The API permits more concurrent requests than the population size of 50, allowing all candidate-generation requests for a generation to be submitted at once; in our experiments, these requests were handled using up to 50 worker threads. \Cref{tab:runtime_breakdown} reports the mean runtime per search run over the ten LLM-guided runs reported in this subsection.

\begin{table}[htbp]
\centering
\caption{Mean wall-clock time per LLM-guided search run over ten runs.}
\label{tab:runtime_breakdown}
\small
\setlength{\tabcolsep}{8pt}
\renewcommand{\arraystretch}{1.08}
\begin{tabular}{@{}lrr@{}}
\toprule
Component & Mean time per run (s) & Share (\%) \\
\midrule
LLM generation and revision & 738.6 & 77.8 \\
Local candidate-rule evaluation & 210.6 & 22.2 \\
\midrule
Total accounted time & 949.2 & 100.0 \\
\bottomrule
\end{tabular}
\end{table}

As shown in \Cref{tab:runtime_breakdown}, LLM generation and revision account for 77.8\% of the wall-clock time, while local candidate-rule evaluation accounts for 22.2\%. Thus, even when all 480 J30 instances are used for training, local candidate-rule evaluation requires less than one quarter of the total runtime, with most of the wall-clock time spent on LLM generation and revision.

The wall-clock time here cannot be compared directly with the original GP implementation of \citet{luo2022efficient}, because the GP library used in that study did not allow fitness evaluations to be executed in parallel. The original configuration required approximately 500 seconds per generation for 50 generations, or about \(500\times50=25{,}000\) seconds per run. Nevertheless, the observed local candidate-rule evaluation time can be used to estimate the time required to evaluate approximately 50,000 candidates under the parallel-computing conditions used in this study. Scaling the observed local candidate-rule evaluation time of 210.6 seconds from a nominal new-candidate evaluation budget of 1,058 to the approximately 50,000-candidate budget of the original GP configuration gives an estimate of \(210.6/1{,}058\times50{,}000=9{,}954\) seconds on the complete J30 set. This estimate includes only candidate-rule evaluation and excludes GP selection, crossover, mutation, and other population operations.

The GP configurations instead use the 96-instance J30 training subset described above, which contains one-fifth of the 480-instance J30 set and was shown by \citet{luo2022efficient} to provide test performance comparable to training on the complete set. Assuming that local candidate-rule evaluation time scales approximately linearly with the number of training instances, the corresponding budget-normalized estimate for the smaller training set is approximately \(9{,}954/5=1{,}991\) seconds.

Overall, the LLM-guided framework reaches the performance range of the large-budget GP search with a nominal new-candidate evaluation budget of approximately 2\% of that used by the GP search. The runtime results also show that local candidate-rule evaluation can be handled efficiently using parallel computing. API usage and monetary cost are reported in \Cref{sec:results_backbones}.

\subsection{Ablation study}
\label{sec:results_ablation}

To assess how the search design contributes to performance, we compare the full framework with five ablation variants. The study focuses on operationally meaningful alternatives rather than every possible combination of the search components. Several components require information or rules from earlier generations. Without cross-generation dependence, for example, there is no earlier rule to pass as a base rule and no earlier fitness information to provide as feedback; Hall-of-Fame references and survivor retention likewise require previously evaluated rules. Testing every possible combination would therefore create configurations that cannot operate as specified or that reduce to the same search procedure. The five operationally meaningful alternatives examined in this experiment are defined below. All variants use the sampling-temperature settings in \Cref{tab:revision_settings}, except where noted.

\begin{description}
\setlength{\itemsep}{0.25em}

\item[\textbf{A1.}] \textbf{Operator-diverse independent sampling.}
Candidates are generated independently from scratch using four separate prompts that request balanced, compact, alternative-driver, and state-gated rules, respectively. No information is passed across generations, and the best sampled rule is returned.

\item[\textbf{A2.}] \textbf{Generic revision prompt.}
The four operator-specific prompts are replaced by one generic prompt that asks the LLM to improve the supplied base rule using the available feedback. The other search components are retained.

\item[\textbf{A3.}] \textbf{Score-only Hall of Fame and reference selection.}
Hall-of-Fame updating is changed to score-only selection. The best-performing members of the resulting score-only archive are supplied as reference rules. The other search components are retained.

\item[\textbf{A4.}] \textbf{No survivor retention.}
Previously evaluated rules are not carried into the next population. Each new population is therefore filled entirely with newly generated candidates.

\item[\textbf{A5.}] \textbf{Plain independent sampling.}
Candidates are generated independently from scratch using one generic prompt and a fixed sampling temperature of 1.0. The best sampled rule is returned.

\end{description}

\Cref{tab:ablation_components} summarizes the search components retained in the full framework and the five variants.

\begin{table}[htbp]
\centering
\begin{threeparttable}
\caption{Search components retained in the full framework and ablation variants.}
\label{tab:ablation_components}
\scriptsize
\setlength{\tabcolsep}{3.2pt}
\renewcommand{\arraystretch}{1.15}
\begin{tabular}{@{}lccccccc@{}}
\toprule
Configuration & \makecell{Cross-generation\\dependence} & \makecell{Base-rule\\passing} & \makecell{Performance\\feedback} & \makecell{Operator-specific\\prompts} & \makecell{HoF\\references} & \makecell{Diversity-aware\\HoF management} & \makecell{Survivor\\retention} \\
\midrule
Full & \checkmark & \checkmark & \checkmark & \checkmark & \checkmark & \checkmark & \checkmark \\
A1 & -- & -- & -- & \checkmark\tnote{a} & -- & -- & -- \\
A2 & \checkmark & \checkmark & \checkmark & -- & \checkmark & \checkmark & \checkmark \\
A3 & \checkmark & \checkmark & \checkmark & \checkmark & \checkmark & -- & \checkmark \\
A4 & \checkmark & \checkmark & \checkmark & \checkmark & \checkmark & \checkmark & -- \\
A5 & -- & -- & -- & -- & -- & -- & -- \\
\bottomrule
\end{tabular}
\begin{tablenotes}[flushleft]
\scriptsize
\item[a] In A1, this check mark denotes four separate from-scratch prompts aligned with the four revision directions. A1 does not use the revision prompts of the full framework.
\end{tablenotes}
\end{threeparttable}
\end{table}

A1 and A5, the two independent-sampling references, each use a fixed budget of 1,058 LLM calls, matching the standard-case count for the full framework when eight survivors are retained in each generation. A2 and A3 retain the same population size and number of generations, but additional calls are required whenever fewer than eight distinct survivors are available. Their mean per-run call counts over the ten runs are 1,152 and 1,154, respectively. A4 retains no survivors and therefore uses 1,250 LLM calls per run.

All configurations use J30 as the training set and are evaluated over ten independent search runs. The full framework serves as the comparison baseline. Each ablation variant is compared with it using a one-sided Mann--Whitney \(U\) test, with the alternative hypothesis that the full framework yields lower AvgDevLB. \Cref{tab:ablation_comparison} reports the test-set performance and LLM-call budget of each configuration. Boldface identifies the lowest mean AvgDevLB in each column.

\begin{table}[htbp]
\centering
\caption{Test-set performance of the full framework and ablation variants.}
\label{tab:ablation_comparison}
\scriptsize
\setlength{\tabcolsep}{6pt}
\renewcommand{\arraystretch}{1.08}
\sisetup{detect-weight=true,mode=text,table-number-alignment=center}
\begin{tabular}{@{}l*{5}{S[table-format=2.3]}c@{}}
\toprule
Configuration & {J60} & {J90} & {J120} & {RG30} & {RG300} & \makecell{LLM calls\\per run} \\
\midrule
Full & \bfseries 4.232 & \bfseries 4.599 & \bfseries 11.743 & \bfseries 6.392 & \bfseries 12.203 & 1,058 \\
A1 & 4.321\sig & 4.678 & 11.825 & 6.818\sig & 13.459\sig & 1,058 \\
A2 & 4.314\sig & 4.676\sig & 11.961\sig & 6.889\sig & 13.606\sig & 1,152 \\
A3 & 4.392\sig & 4.737\sig & 12.077\sig & 7.065\sig & 14.119\sig & 1,154 \\
A4 & 4.385\sig & 4.730\sig & 12.035\sig & 6.995\sig & 13.510\sig & 1,250 \\
A5 & 4.300 & 4.704\sig & 11.782 & 6.946\sig & 13.197\sig & 1,058 \\
\bottomrule
\end{tabular}
\end{table}

\Cref{tab:ablation_comparison} shows a consistent pattern across the five test sets: every ablation variant yields a higher mean AvgDevLB than the full framework. The performance losses are generally larger on RG30 and RG300 than on the PSPLIB test sets, suggesting that the advantage of the full framework becomes more evident when the designed rules are applied to a different instance family.

Among the variants that retain the cross-generation search structure, A3 shows the largest deterioration on every test set. A2 and A4 are also significantly worse than the full framework on all five test sets. These results favor the operator-specific prompts, diversity-aware Hall-of-Fame management, and survivor retention used in the full framework.

These losses cannot be attributed to smaller call budgets. A2, A3, and A4 use 1,152, 1,154, and 1,250 LLM calls per run, respectively, all exceeding the standard-case count of 1,058 for the full framework. A3 and A4 show particularly large performance losses despite also using the two largest call budgets.

A1 and A5 remain relatively close to the full framework on the PSPLIB test sets but deteriorate more clearly on RG30 and RG300. Notably, A5 does not produce the largest performance loss despite retaining none of the search components listed in \Cref{tab:ablation_components}. The effects of the search components therefore depend on the overall search configuration rather than accumulating independently. Because some variants change more than one component at the same time, their results should be interpreted at the variant level and cannot be assigned separately to individual components.

\subsection{Expanded training and comparison with the best-performing GP rule}
\label{sec:results_expanded}

To evaluate the effect of expanding the training set, we compare J30-only training with J30+RG30 training. RG30 contains projects generated using RanGen and spans broader ranges of network-topology indicators than J30, thereby exposing the search to a broader range of project structures. Because RG30 is included in the expanded training set, the comparison focuses on J60, J90, J120, and RG300, which remain outside the training data under both settings. Each training setting is evaluated over ten independent runs.

We also compare the expanded-training results with GP1, the best-performing GP-designed priority rule reported by \citet{luo2022efficient}. GP1 was obtained using the 186-instance J30+RG30 training subset adopted in that study, which combines 96 J30 instances, comprising the first two instances from each of the 48 PSPLIB parameter combinations, with 90 instances sampled at intervals of ten from RG30 Set~1. As discussed in \Cref{sec:results_baselines}, \citet{luo2022efficient} showed that carefully sampled training sets could retain test-set performance comparable to that obtained from much larger training sets. Using these smaller training sets also substantially reduced the number of schedule calculations required during GP evolution.

\Cref{tab:expanded_training_gp1} combines three comparisons. The first two columns report the ten-run means for J30 and J30+RG30 training. The next four columns reproduce the traditional-rule results from \Cref{tab:baseline_comparison}, and the final block compares one selected LLM rule with GP1. For the expanded-training comparison, J30 training is the baseline for one-sided Mann--Whitney \(U\) tests with the alternative hypothesis that J30+RG30 training yields lower AvgDevLB. For the selected-rule comparison, the LLM rule is chosen retrospectively as the run with the lowest mean AvgDevLB across J60, J90, J120, and RG300. The W/T/L column reports the numbers of instances on which this rule obtains a lower, equal, or higher makespan than GP1. These counts are descriptive. Boldface marks the lower AvgDevLB in the selected-rule comparison.

\begin{table}[htbp]
\centering
\caption{Effect of expanded training and comparison with traditional priority rules and GP1.}
\label{tab:expanded_training_gp1}
\scriptsize
\setlength{\tabcolsep}{3.2pt}
\renewcommand{\arraystretch}{1.08}
\sisetup{detect-weight=true,mode=text,table-number-alignment=center}

\begin{tabular*}{\textwidth}{
@{\extracolsep{\fill}}
l
*{8}{S[table-format=2.3]}
c
@{}
}
\toprule
& \multicolumn{2}{c}{Ten-run mean}
& \multicolumn{4}{c}{Traditional references}
& \multicolumn{3}{c}{Selected-rule comparison} \\
\cmidrule(lr){2-3}
\cmidrule(lr){4-7}
\cmidrule(lr){8-10}

Test set
& {J30}
& {J30+RG30}
& {WCS}
& {LFT}
& {SinglePR}
& {AllPR}
& {\makecell{Selected\\LLM rule}}
& {GP1}
& {W/T/L} \\
\midrule

J60
& 4.232\sig
& 4.128
& 5.746
& 4.555
& 4.372
& 3.334
& \bfseries 4.081
& 4.227
& 80/342/58 \\

J90
& 4.599\sig
& 4.479
& 5.742
& 4.947
& 4.845
& 4.038
& 4.477
& \bfseries 4.368
& 68/336/76 \\

J120
& 11.743\sig
& 11.304
& 13.203
& 12.221
& 12.166
& 10.604
& \bfseries 11.304
& 11.391
& 206/200/194 \\

RG300
& 12.203\sig
& 11.645
& 13.656
& 13.170
& 13.012
& 11.748
& \bfseries 11.476
& 11.843
& 249/96/135 \\

\bottomrule
\end{tabular*}
\end{table}

The expanded-training comparison in \Cref{tab:expanded_training_gp1} shows that adding RG30 to the training set significantly lowers the mean AvgDevLB on all four test sets. The largest improvement occurs on RG300, consistent with the inclusion of RanGen-generated projects in the training set, while the reductions on J60, J90, and J120 show that the benefit also extends to the PSPLIB test sets.

The selected LLM rule obtains lower AvgDevLB than WCS, LFT, and SinglePR on all four test sets. AllPR obtains lower AvgDevLB on J60, J90, and J120, but this ordering reverses on RG300: the selected LLM rule obtains an AvgDevLB of 11.476\%, compared with 11.748\% for AllPR, despite using only a single priority rule. Relative to GP1, the ten-run mean is lower on J60, J120, and RG300 and slightly higher on J90. The result is consistent across the independently generated LLM rules: all ten outperform GP1 on J60 and J120, nine do so on RG300, and none does so on J90. The advantage on the three datasets therefore extends beyond the selected best rule.

The selected-rule comparison follows the same dataset-level pattern. The selected LLM rule outperforms GP1 on J60, J120, and RG300, while GP1 retains a small advantage on J90. The clearest instance-level advantage occurs on RG300, with 249 wins and 135 losses. The PSPLIB datasets contain substantially more ties; on J90, for example, the selected rule records 68 wins and 76 losses, indicating only a small instance-level difference. The selected expanded-training rule is further evaluated on the large-project datasets in \Cref{sec:results_large_projects}, without retraining or a new rule-design search.

\subsection{Backbone robustness and API usage}
\label{sec:results_backbones}

This subsection examines robustness across LLM backbones and reports the associated API usage and cost. All experiments use J30+RG30 training, with only the backbone varied. Gemini Flash, Claude Haiku 4.5, GPT-5.4-mini, and DeepSeek V4-Flash are each evaluated over ten independent runs, with Gemini Flash serving as the reference for two-sided Mann--Whitney \(U\) tests. GPT-5.5 and DeepSeek V4-Pro are included as higher-cost premium LLMs and are each evaluated once as exploratory references. \Cref{tab:backbone_performance} reports the test-set results. Boldface identifies the lowest mean AvgDevLB among the four ten-run backbones for each test set; the exploratory single-run results are excluded from this ranking.

\begin{table}[htbp]
\centering
\caption{Test-set performance across LLM backbones.}
\label{tab:backbone_performance}
\scriptsize
\setlength{\tabcolsep}{6pt}
\renewcommand{\arraystretch}{1.08}
\sisetup{detect-weight=true,mode=text,table-number-alignment=center}
\begin{tabular}{lc*{4}{S[table-format=2.3]}}
\toprule
Backbone & Runs & {J60} & {J90} & {J120} & {RG300} \\
\midrule
Gemini Flash & 10 & 4.128 & 4.479 & 11.304 & 11.645 \\
Claude Haiku 4.5 & 10 & 4.139 & 4.455 & \bfseries 11.268 & 11.827 \\
GPT-5.4-mini & 10 & 4.202 & 4.556 & 11.654\sig & 12.654\sig \\
DeepSeek V4-Flash & 10 & \bfseries 4.123 & \bfseries 4.409\sig & 11.270 & \bfseries 11.423\sig \\
\midrule
GPT-5.5 & 1 & 4.099 & 4.492 & 11.291 & 11.531 \\
DeepSeek V4-Pro & 1 & 4.239 & 4.439 & 11.133 & 11.420 \\
\bottomrule
\end{tabular}
\end{table}

The ten-run results show that the framework is not tied to a single backbone, although rule quality varies across models. DeepSeek V4-Flash obtains the lowest mean AvgDevLB on J60, J90, and RG300, with significant improvements over Gemini Flash on J90 and RG300. Claude Haiku 4.5 performs best on J120 and remains close to Gemini Flash elsewhere, whereas GPT-5.4-mini performs significantly worse on J120 and RG300.

The two exploratory premium-model runs are competitive with the ten-run backbone results but do not show a consistent performance advantage. DeepSeek V4-Pro obtains the lowest AvgDevLB on J120 and RG300 among the reported backbone results, while GPT-5.5 obtains the lowest value on J60; neither performs best across all four test sets.

\Cref{tab:backbone_usage} reports input-token usage, provider-reported output-token usage, and approximate API cost per search run. For the four main backbones, the reported values are averages over ten independent runs; for GPT-5.5 and DeepSeek V4-Pro, they correspond to the single exploratory run. Because the providers differ in token accounting and billing conventions, the usage statistics are reported as provided by the corresponding APIs.

\begin{table}[htbp]
\centering
\caption{API usage and cost per search run across LLM backbones.}
\label{tab:backbone_usage}
\scriptsize
\setlength{\tabcolsep}{5pt}
\renewcommand{\arraystretch}{1.08}
\sisetup{table-number-alignment=center}
\begin{tabular*}{\textwidth}{@{\extracolsep{\fill}}lcS[table-format=1.2]S[table-format=1.2]c@{}}
\toprule
Backbone & Runs & {\makecell{Input tokens\\(millions)}} & {\makecell{Provider-reported output\\tokens (millions)}} & {\makecell{Approx.\\cost}} \\
\midrule
Gemini Flash & 10 & 4.36 & 2.67 & USD 10.16 \\
Claude Haiku 4.5 & 10 & 4.75 & 2.92 & USD 19.35 \\
GPT-5.4-mini & 10 & 4.59 & 1.81 & USD 10.36 \\
DeepSeek V4-Flash & 10 & 5.24 & 1.79 & CNY 6.87 \\
\midrule
GPT-5.5 & 1 & 4.18 & 0.91 & USD 44.54 \\
DeepSeek V4-Pro & 1 & 5.53 & 2.46 & CNY 25.36 \\
\bottomrule
\end{tabular*}
\end{table}

Among the four ten-run backbones, mean input-token usage per run is similar, while API cost varies more substantially across providers. Gemini Flash and GPT-5.4-mini have similar USD costs, Claude Haiku 4.5 is more expensive, and DeepSeek V4-Flash remains inexpensive under its provider's CNY pricing. The two exploratory premium LLMs are considerably more expensive per run than their lower-cost counterparts without showing a consistent performance advantage. Because providers account for cached inputs, reasoning tokens, and other usage categories differently, these figures should not be interpreted as a precise cross-provider cost ranking.

\subsection{Performance on large projects}
\label{sec:results_large_projects}

We further evaluate whether rules designed using 30-activity projects retain their effectiveness when applied to substantially larger projects. The analysis uses the Multi-projects, LPP, and LPSP datasets described in \Cref{sec:benchmark_design}, each of which contains projects with more than 1,000 activities. Following the evaluation protocol established for these datasets, performance is reported using AvgDevBKS, with lower values indicating better schedules.

Two LLM-generated rules are included in the comparison. The Gemini Flash rule is the same rule used in the head-to-head comparison with GP1 in \Cref{sec:results_expanded}; the DeepSeek V4-Flash rule was selected from its ten J30+RG30 runs using the same four-test-set mean criterion. Both selections preceded the large-project evaluation.

\Cref{tab:large_project_performance} compares the two selected LLM rules with three configurations of the decomposition-based GA proposed by \citet{debels2007decomposition}, GP1, the traditional priority-rule baselines, and the best regression-based heuristic reported by \citet{luo2026automated} for each dataset. The GA configurations are denoted by DV$_{1000}$, DV$_{5000}$, and DV$_{50000}$, where the subscript gives the number of schedule evaluations per instance. Once a rule or fitted scoring model has been fixed, repeated application of the same non-GA method to a given instance produces the same schedule under the common schedule-generation and tie-breaking protocol. Because the GA produces stochastic run-to-run variation, its entries report mean AvgDevBKS over ten independent runs. Boldface identifies the best non-GA result for each dataset.

\begin{table}[htbp]
\centering
\caption{Performance on large-project datasets (AvgDevBKS, \%).}
\label{tab:large_project_performance}
\scriptsize
\setlength{\tabcolsep}{2.4pt}
\renewcommand{\arraystretch}{1.10}
\sisetup{detect-weight=true,mode=text,table-number-alignment=center}
\begin{tabular*}{\textwidth}{@{\extracolsep{\fill}}l*{11}{S[table-format=1.2]}@{}}
\toprule
& \multicolumn{3}{c}{Genetic algorithm} & \multicolumn{1}{c}{GP} & \multicolumn{4}{c}{Priority-rule baselines} & \multicolumn{1}{c}{Regression} & \multicolumn{2}{c}{Selected LLM rules} \\
\cmidrule(lr){2-4}\cmidrule(lr){5-5}\cmidrule(lr){6-9}\cmidrule(lr){10-10}\cmidrule(lr){11-12}
Dataset & {DV$_{1000}$} & {DV$_{5000}$} & {DV$_{50000}$} & {GP1} & {WCS} & {LFT} & {SinglePR} & {AllPR} & {\makecell{Best\\regression}} & {\makecell{Gemini\\Flash}} & {\makecell{DeepSeek\\V4-Flash}} \\
\midrule
Multi-projects & 1.72 & 1.16 & 0.11 & 1.00 & 1.27 & 1.45 & 1.18 & 0.97 & 1.02 & 0.81 & \bfseries 0.72 \\
LPP & 3.75 & 2.78 & 1.33 & 2.87 & 4.02 & 3.60 & 3.53 & 2.98 & 3.17 & 2.54 & \bfseries 2.33 \\
LPSP & 1.31 & 0.81 & 0.14 & 3.67 & 2.78 & 2.37 & 2.27 & 2.16 & 2.17 & \bfseries 2.06 & 2.08 \\
\bottomrule
\end{tabular*}
\end{table}

Turning first to the non-GA results, both selected LLM rules achieve lower AvgDevBKS than GP1, all single-rule traditional baselines, and the best regression-based heuristic on each dataset. DeepSeek V4-Flash achieves the best non-GA result on Multi-projects and LPP, whereas Gemini Flash does so on LPSP.

The comparison with AllPR reveals a notable reversal. AllPR outperforms the J30-trained LLM rules on J60, J90, J120, RG30, and RG300. After expanding the training set, AllPR remains better on J60, J90, and J120, but the ordering reverses on RG300. On all three large-project datasets, both selected LLM rules also obtain lower AvgDevBKS than AllPR. Thus, the advantage of applying all 36 traditional rules and retaining the best schedule does not persist on the large-project datasets.

Compared with the GA references, both selected LLM rules obtain lower AvgDevBKS than DV$_{1000}$ and DV$_{5000}$ on Multi-projects and LPP. Both datasets are characterized by highly parallel project networks. Such structures impose relatively weak precedence restrictions and admit more precedence-feasible activity orderings, making high-quality solutions harder to identify under limited evaluation budgets \citep{coelho2020going}. DV$_{50000}$ performs better on both datasets. Thus, on these two highly parallel datasets, increasing the per-instance search to 5,000 schedule evaluations is still insufficient for the GA to match the selected LLM rules, whereas 50,000 evaluations reverse the ordering. On LPSP, all three GA configurations outperform the selected LLM rules.

Overall, the selected rules remain effective far beyond the 30-activity projects used for their design, although their relative performance against per-instance GA search depends on project structure and search budget.

\section{Analysis and discussion of learned priority rules}
\label{sec:rule_analysis}

This section goes beyond aggregate performance to examine the structure and decision behavior of the priority rules and to identify the scheduling logic represented by the generated rules. The analysis covers 71 priority rules: 52 LLM-generated rules, 10 GP-designed rules, and 9 traditional single priority rules. Twenty of the LLM-generated rules were produced with Gemini Flash, split evenly between J30 and J30+RG30 training. The remaining 32 come from the backbone experiments: 10 each from Claude Haiku 4.5, GPT-5.4-mini, and DeepSeek V4-Flash, together with the two exploratory single-run rules generated with the premium LLMs GPT-5.5 and DeepSeek V4-Pro. \Cref{sec:rule_structure,sec:decision_behavior} present the structural and behavioral analyses, while \Cref{sec:deployment_implications} discusses their implications for rule design and deployment.

\subsection{Rule structure and contribution of project-state information}
\label{sec:rule_structure}

This subsection compares the structure of the rules generated under J30 and J30+RG30 training and examines whether project-level indicators and progress-dependent switching contribute to rule performance. The structural comparison uses three equally sized groups: the 10 Gemini Flash rules trained on J30, the 10 Gemini Flash rules trained on J30+RG30, and the 10 GP-designed rules produced by the original GP configuration using the 96-instance J30 training subset. Comparing the two Gemini Flash groups keeps the LLM backbone and search framework unchanged and varies only the training setting. The GP group provides a benchmark generated by a different rule-design method.

Rules generated with the alternative LLM backbones, together with the traditional rules, are included in the behavioral analysis in \Cref{sec:decision_behavior}. \Cref{tab:rule_taxonomy} reports the prevalence of progress-dependent switching, the mean number of project-level indicators, the use of sine or cosine, and the mean number of activity attributes in each group. Progress-dependent switching and project-level indicators do not apply to GP because these inputs were absent from its terminal set.

\begin{table}[htbp]
\centering
\caption{Structural taxonomy of the three ten-run rule groups.}
\label{tab:rule_taxonomy}
\scriptsize
\setlength{\tabcolsep}{5pt}
\renewcommand{\arraystretch}{1.10}
\begin{tabular*}{\textwidth}{@{\extracolsep{\fill}}lcccc@{}}
\toprule
Rule group
& \makecell{Rules with\\progress switching}
& \makecell{Project indicators\\per rule (mean)}
& \makecell{Rules using\\sine or cosine}
& \makecell{Activity attributes\\per rule (mean)} \\
\midrule
LLM, J30+RG30 training & 9/10 & 3.9 & 0/10 & 6.8 \\
LLM, J30 training & 2/10 & 1.7 & 0/10 & 6.3 \\
GP, original configuration & N/A & N/A & 8/10 & 6.5 \\
\bottomrule
\end{tabular*}
\end{table}

As shown in \Cref{tab:rule_taxonomy}, the three groups use similar numbers of activity attributes per rule. The performance improvement under expanded training is therefore not accompanied by a substantial increase in the number of activity attributes used. Among the J30+RG30-trained LLM rules, all ten use latest finish time, latest start time, resource demand, and total successor processing time; nine also use cumulative resource work content and the length of the longest path following the activity. These inputs span the scheduling, resource, and network categories defined in \Cref{tab:activity_inputs}, representing timing urgency, resource pressure, and downstream precedence information, respectively. Taken as a group, the ten GP rules use 21 distinct activity attributes, whereas the ten J30+RG30-trained LLM rules use nine. Eight of the ten GP rules also contain sine or cosine, while none of the LLM rules do. The LLM rules therefore rely more consistently on a smaller shared set of activity attributes, while the GP rules vary more in both attribute selection and functional form.

Among the LLM rules, the J30+RG30-trained rules use project-level indicators and progress-dependent switching more frequently than the J30-trained rules. Nine of the ten J30+RG30-trained rules use progress-dependent switching between two phases, compared with two of the ten J30-trained rules, and the mean number of project-level indicators per rule rises from 1.7 to 3.9. These indicators generally adjust the weights assigned to activity-level terms according to project characteristics, rather than serving as independent priority measures. In rules with progress-dependent switching, the expression used before the progress threshold typically places greater emphasis on latest finish time, resource demand, and the length of the longest path following the activity. The expression used after the threshold places greater emphasis on latest start time and cumulative resource work content.

To illustrate how these activity-level and project-state components are combined in an executable rule, we examine a formula that was rediscovered independently across runs and backbones. Two independent J30+RG30 runs produced the same formula, and the GPT-5.5 run produced an algebraically equivalent version using different variable names. Following the notation introduced in \Cref{sec:rule_representation}, lower-case roman symbols denote the scaled rule inputs supplied to the priority function. For an eligible activity \(i\), define
\begin{align}
\mathrm{impact}_i
&=
\mathrm{mtspt}_i
+
\mathrm{lpf}_i(1.1+\mathrm{sp}),
\label{eq:reference_rule_impact}\\
\mathrm{pressure}_i
&=
\left(\mathrm{grd}_i+0.05\,\mathrm{crwc}_i\right)
(\mathrm{rc}+\mathrm{rf})(1.1-\mathrm{rs}).
\label{eq:reference_rule_pressure}
\end{align}

The priority score at decision point \(t\) is
\begin{equation}
q_{it}
=
\begin{cases}
\mathrm{lf}_i-\mathrm{impact}_i-\mathrm{pressure}_i,
& \mathrm{progress}_t<0.75,\\
\mathrm{ls}_i-\mathrm{impact}_i-0.5\,\mathrm{pressure}_i,
& \mathrm{progress}_t\geq 0.75.
\end{cases}
\label{eq:reference_rule}
\end{equation}

Because lower scores receive higher priority, the \(\mathrm{impact}_i\) term favors activities with larger \(\mathrm{mtspt}_i\) and \(\mathrm{lpf}_i\) values. The project-network indicator \(\mathrm{sp}\) increases the contribution of \(\mathrm{lpf}_i\) in more serial project networks. The \(\mathrm{pressure}_i\) term favors activities with greater duration-weighted resource demand and greater resource work content among their immediate successors, while \(\mathrm{rc}\), \(\mathrm{rf}\), and \(\mathrm{rs}\) adjust this contribution according to the resource characteristics of the project. When \(\mathrm{progress}_t<0.75\), the rule uses latest finish time and the full pressure term. After this threshold is reached, it switches to latest start time and halves the contribution of resource pressure.

The taxonomy shows that project-level indicators and progress-dependent switching occur frequently in the J30+RG30-trained rules. Their presence in the formulas, however, does not establish whether they contribute to performance. We therefore construct two counterfactual versions while leaving all activity-level terms and coefficients unchanged. First, the project-level indicators are replaced by their empirical means on J30, while progress remains dynamic. Second, progress is also fixed at its J30 empirical mean, so that the progress-based rules use the same branch throughout schedule construction. The second comparison therefore measures the additional effect of removing progress-dependent switching after the project indicators have already been held constant and is applied to the nine rules containing this structure. Performance is measured by mean AvgDevLB across J60, J90, J120, and RG300. \Cref{tab:state_counterfactual} summarizes the results.

\begin{table}[htbp]
\centering
\caption{Counterfactual effects of project-state information on mean AvgDevLB (\%).}
\label{tab:state_counterfactual}
\scriptsize
\setlength{\tabcolsep}{10pt}
\renewcommand{\arraystretch}{1.10}
\begin{tabular}{@{}lcccc@{}}
\toprule
& & \multicolumn{2}{c}{Mean AvgDevLB (\%)} & \\
\cmidrule(lr){3-4}
Modification & Rules worsened & Before & After & Increase \\
\midrule
Project indicators held constant & 10/10 & 7.889 & 8.041 & +0.152 \\
Progress also held constant & 9/9 & 8.024 & 8.154 & +0.130 \\
\bottomrule
\end{tabular}
\end{table}

As shown in \Cref{tab:state_counterfactual}, holding the project-level indicators constant increases the mean AvgDevLB of the ten rules from 7.889\% to 8.041\%. Among the nine rules with progress-dependent switching, the corresponding indicator-fixed mean AvgDevLB is 8.024\%; additionally holding progress constant increases it to 8.154\%. Every applicable rule deteriorates under both modifications. These results suggest that the use of project-level indicators and progress-dependent switching contributed to the performance of the J30+RG30-trained rules.

\subsection{Decision behavior and test-set performance}
\label{sec:decision_behavior}

Rules with different formulas may still make similar activity choices, so structural analysis alone does not fully describe the scheduling logic represented by the generated rules. This subsection therefore examines rule behavior in two steps. First, we apply the analyzed rules to the same decision states and examine whether independently generated rules make similar activity choices. Second, we use a repeatedly rediscovered LLM rule as a behavioral reference and examine whether behavioral distance from this reference is related to test-set performance.

\textit{Behavioral similarity across rule groups.} To ensure that all analyzed rules face the same decision states, we construct decision-state pools from schedules generated by fixed anchor rules. At each decision point \(t\) with at least two eligible activities, the recorded decision state is
\begin{equation}
\delta_t=
\left(
\left\{(i,\mathbf{x}_i):i\in E_t\right\},
\mathbf{z}_t
\right),
\label{eq:decision_state}
\end{equation}
which records the eligible activities and their activity-level inputs together with the project-level and schedule-construction inputs at that decision point. For each pool, the states collected from the three anchor-rule runs are combined and duplicates are removed. Every analyzed rule is then applied to exactly the same retained states, and its top-priority activity is recorded. The resulting differences therefore reflect how the rules rank the same eligible activities, rather than differences in the decision states that would arise if each rule generated its own schedule.

Let \(\mathcal{D}\) denote a decision-state pool. We define the behavioral distance between two priority rules on \(\mathcal{D}\) as the proportion of states in which the two rules select different top-priority activities. Two pools are used in the analysis. The J30 pool is generated by applying LFT, MTS, and GRD to all 480 J30 instances and contains 35,824 unique decision states. The full test-set pool is constructed in the same way using all instances in J60, J90, J120, and RG300. Its anchors are the repeatedly rediscovered LLM rule in \Cref{eq:reference_rule}, the best-performing of the ten GP-designed rules included in the analysis, and LFT, representing the LLM, GP, and traditional-rule groups, respectively. The resulting pool contains 857,006 unique decision states. The repeatedly rediscovered LLM rule also serves as the behavioral reference rule in the subsequent distance--performance analysis. Algebraically equivalent versions of this rule were generated independently across runs and LLM backbones, indicating that it captures a recurring decision logic. It was therefore selected as the reference independently of test-set performance.

\begin{figure}[!t]
\centering
\includegraphics[width=1\textwidth]{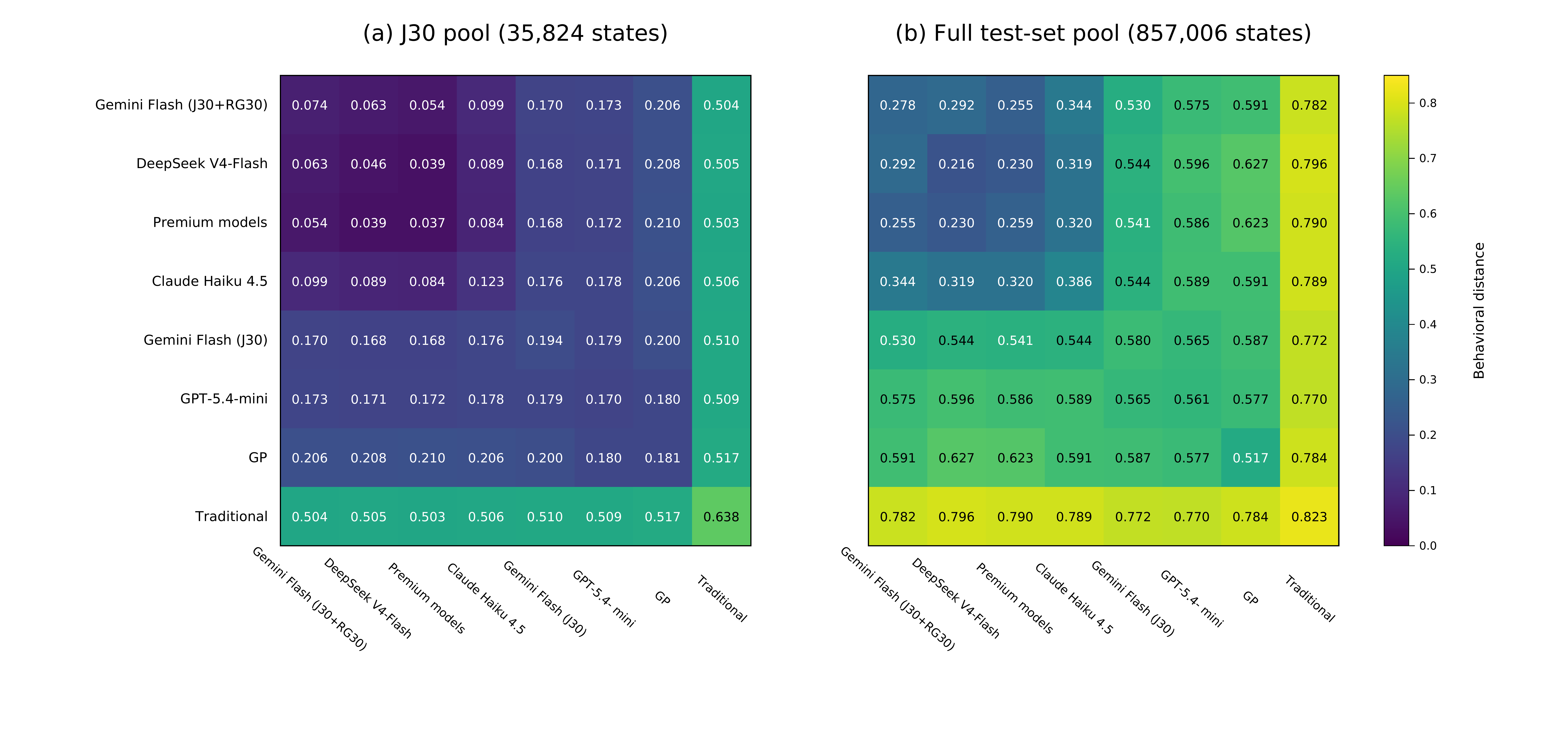}
\caption{Behavioral distance within and across rule groups on the J30 and full test-set pools.}
\label{fig:group_behavioral_distance}
\end{figure}

We compare within-group and between-group behavioral distances to assess rule similarity within and across groups. Repeating the comparison on both decision-state pools shows whether these relative similarities change with the decision states used. \Cref{fig:group_behavioral_distance} presents the two matrices using a common group order and color scale. Off-diagonal cells average the distances over all rule pairs drawn from the corresponding groups, whereas diagonal cells average the distances over all distinct rule pairs within the same group; self-comparisons are excluded.

As shown in \Cref{fig:group_behavioral_distance}, absolute behavioral distances are generally higher on the full test-set pool than on the J30 pool, but the relative relationships among the rule groups are similar. Among the four backbone groups evaluated over ten runs under J30+RG30 training, the ordering of within-group distance is identical in both pools: DeepSeek V4-Flash has the lowest distance, followed by Gemini Flash, Claude Haiku 4.5, and GPT-5.4-mini. In the between-group comparisons, Gemini Flash (J30) is farther from the other LLM groups than Gemini Flash (J30+RG30) is, but it remains closer to every LLM group than to either GP or the traditional rules in both pools.

\textit{Relationship between behavior and performance.} The primary analysis uses behavioral distances measured on the J30 pool. Its decision states come from J30, which is separate from the four test sets used to measure performance, and the behavioral reference rule is not one of its anchors. We repeat the analysis using the full test-set pool as a consistency check. In this second pool, the decision states come from the same four test sets, and the behavioral reference rule is one of the anchors.

The correlation analysis covers the 62 automatically designed rules; the nine traditional rules are retained only in the group-level comparison. For each rule, test-set performance is summarized by the unweighted mean AvgDevLB across J60, J90, J120, and RG300. We also report separate correlations for two training subsets to examine whether the association persists within each training setting. The J30-trained subset contains the 10 LLM rules and the 10 GP rules from the original configuration reported in \Cref{sec:results_baselines}. The J30+RG30-trained subset contains the 10 Gemini Flash rules from \Cref{sec:results_expanded} and the 32 rules from \Cref{sec:results_backbones}. Spearman's rank correlation is used to test whether greater behavioral distance from the reference is associated with worse test-set performance. Because lower AvgDevLB indicates better performance, a positive coefficient means that rules farther from the reference tend to perform worse. \Cref{tab:distance_performance_correlation} reports the results.

\begin{table}[htbp]
\centering
\caption{Spearman correlations between reference-rule distance and test-set performance.}
\label{tab:distance_performance_correlation}
\scriptsize
\setlength{\tabcolsep}{4pt}
\renewcommand{\arraystretch}{1.15}
\begin{tabular*}{0.96\textwidth}{@{\extracolsep{\fill}}lcccc@{}}
\toprule
\makecell[l]{State pool used\\for distance}
& \makecell{All 62\\rules}
& \makecell{Excluding the three\\zero-distance rules}
& \makecell{J30-trained\\rules}
& \makecell{J30+RG30-trained\\rules} \\
\midrule
J30 pool & 0.856 & 0.873 & 0.615 & 0.698 \\
Full test-set pool & 0.872 & 0.892 & 0.671 & 0.709 \\
\bottomrule
\end{tabular*}
\end{table}

As shown in \Cref{tab:distance_performance_correlation}, the primary analysis based on the J30 pool shows a strong positive association between behavioral distance from the reference and mean test-set AvgDevLB. The association remains strong after excluding the three zero-distance rules and remains positive when the J30-trained and J30+RG30-trained rules are analyzed separately.

\begin{figure}[htbp]
\centering
\includegraphics[width=0.7\textwidth]{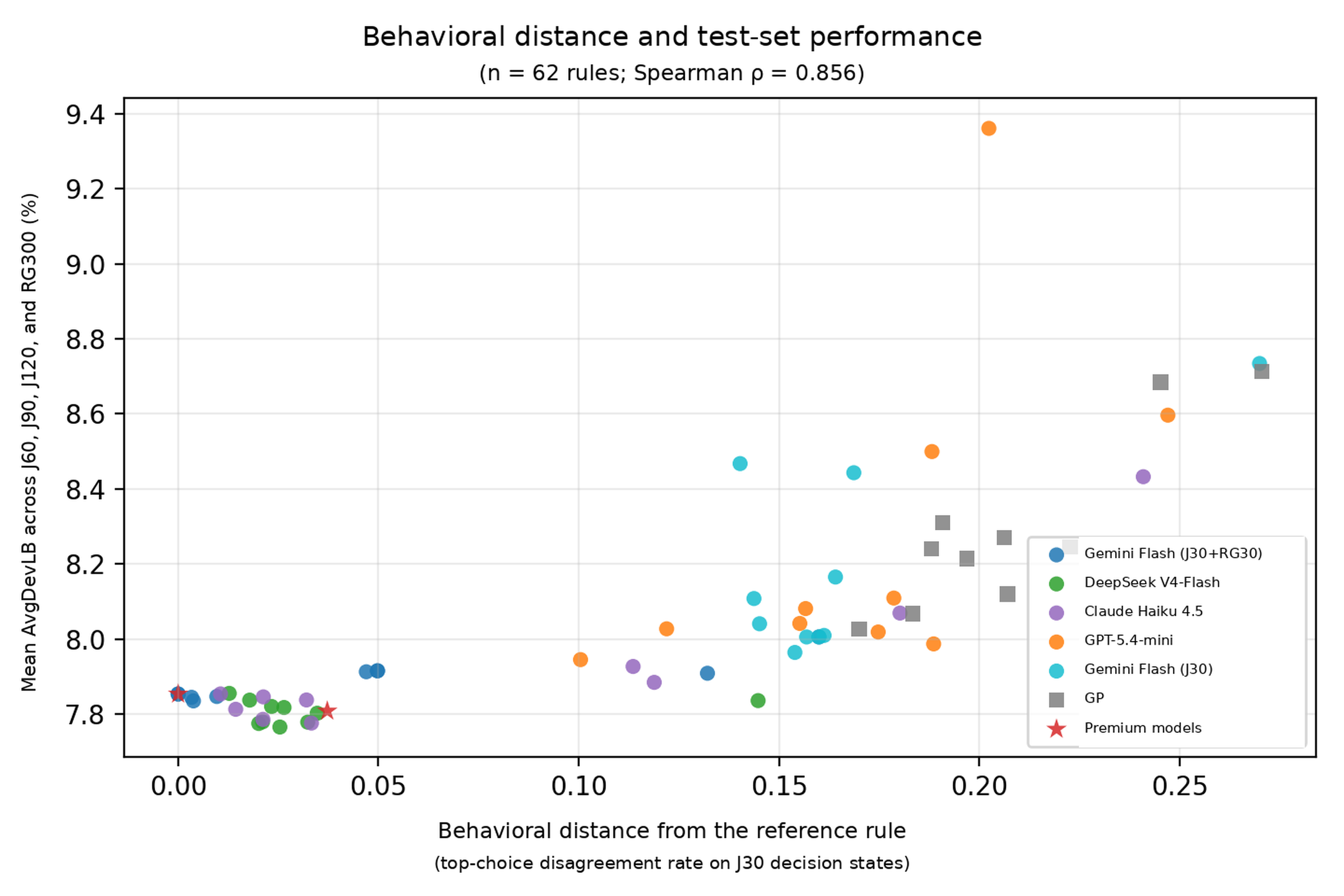}
\caption{Behavioral distance measured on the J30 pool versus mean test-set AvgDevLB.}
\label{fig:distance_performance_overall}
\end{figure}

Using behavioral distances measured on the J30 pool, \Cref{fig:distance_performance_overall} visualizes the relationship for all 62 automatically designed rules. Rules with lower mean test-set AvgDevLB are concentrated closer to the reference behavior, while rules farther from it tend to have higher AvgDevLB.

\begin{figure}[!htbp]
\centering
\includegraphics[width=0.82\textwidth]{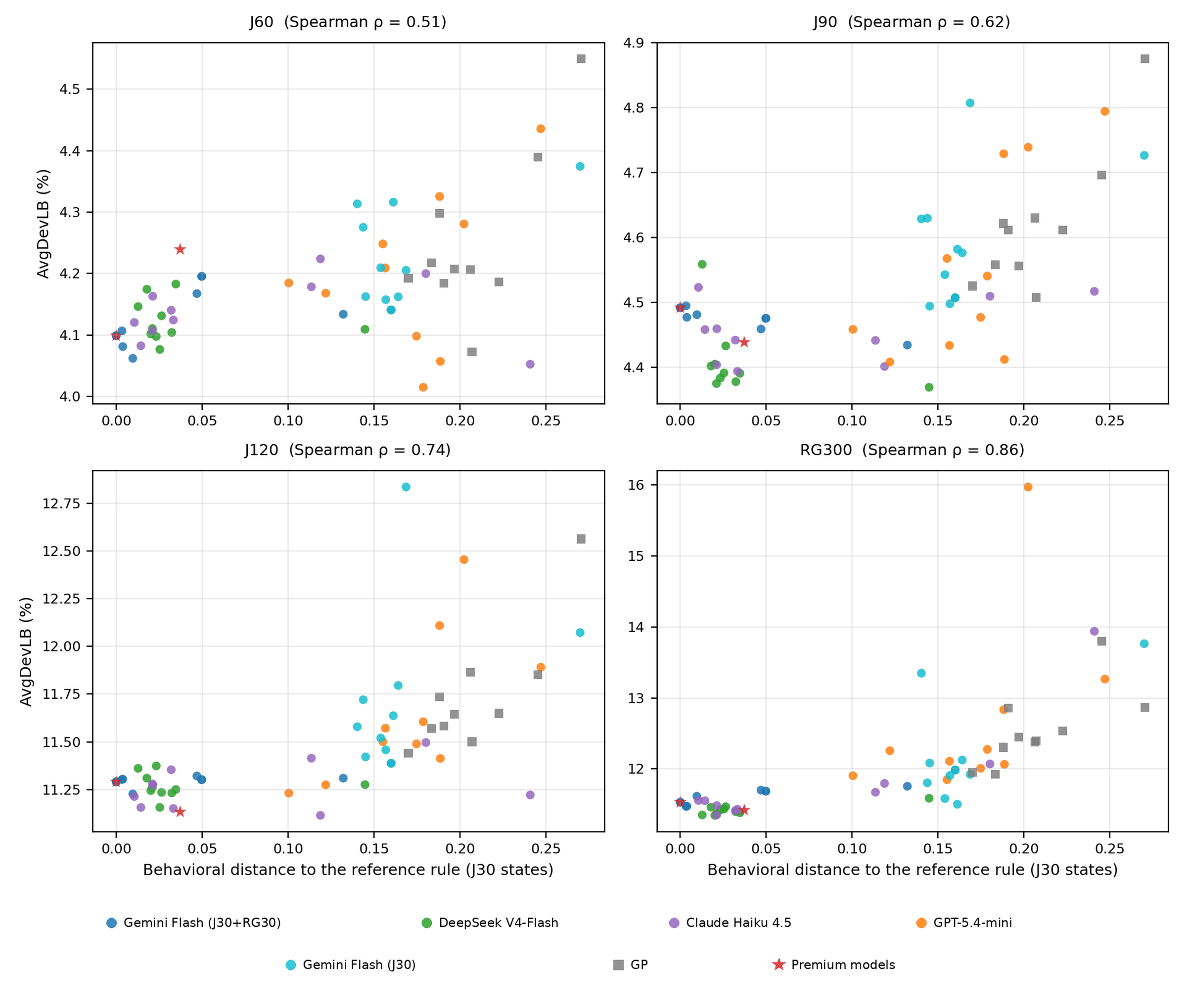}
\caption{Behavioral distance measured on the J30 pool versus AvgDevLB by test set.}
\label{fig:distance_performance_by_dataset}
\end{figure}

\Cref{fig:distance_performance_by_dataset} shows the same relationship separately for J60, J90, J120, and RG300. The association is positive on all four test sets and is strongest on J120 and RG300.

The second row of \Cref{tab:distance_performance_correlation}, based on the full test-set pool, shows a similar relationship. In a separate comparison, the rankings of the 62 rules by behavioral distance under the two pools are highly aligned, with a Spearman correlation of 0.952. The relationship between behavioral proximity and performance is therefore not specific to the J30 decision-state pool. Behavioral distance may provide useful information for rule screening, selection, and diversity management, although its use directly within the search remains a topic for future research.

\subsection{Implications for priority-rule design and deployment}
\label{sec:deployment_implications}

\textit{Rule-design implications.} Among the rules generated in this study, the J30+RG30-trained rules use project-level indicators and progress-dependent switching more frequently than the J30-trained rules. The counterfactual analysis suggests that both elements contributed to their performance. This indicates that allowing a priority rule to respond to project characteristics and schedule-construction progress can improve scheduling performance. A separate implication follows from the repeated discovery of algebraically equivalent rules and the relationship between behavioral similarity and test-set performance. Behavioral similarity may therefore help screen candidates, select final and reference rules, and maintain diversity in future searches.

\textit{Deployment implications.} All LLM-related computational effort is confined to the offline rule-design stage. Once selected, the rule is a fixed deterministic priority function composed of arithmetic operations and conditional expressions and is used directly to rank eligible activities during schedule construction. Scheduling a new project therefore requires no further LLM calls or rule-design search; deployment involves only direct priority-score calculations, as with traditional priority-rule heuristics. This separation is useful for repeated scheduling and rescheduling and becomes especially valuable for large, highly parallel projects.

As discussed in \Cref{sec:results_large_projects}, such project structures admit a much larger set of precedence-feasible activity orderings. As the search space expands, a fixed per-instance evaluation budget provides increasingly limited search coverage, so budgets that are sufficient to identify high-quality schedules in smaller search spaces may no longer be adequate. On the two highly parallel large-project datasets, the selected LLM rules outperform AllPR. They also outperform the GA configurations using 1,000 and 5,000 schedule evaluations per instance, with only the 50,000-evaluation GA obtaining better results. These findings suggest that an offline-designed reusable priority rule can be particularly useful for large, highly parallel projects when obtaining high-quality solutions through per-instance search would require substantially larger evaluation budgets.

\section{Conclusion and future research}
\label{sec:conclusion}

This study examined whether an LLM can serve as the primary search mechanism in an offline process for automatically designing reusable priority rules for the RCPSP. The proposed framework uses schedule-quality feedback to guide LLM-based generation and revision of candidate rules within a population-based search. Each search run returns one deterministic priority function that can subsequently be applied to unseen projects without a new rule-design search or further LLM calls.

The computational evidence shows that the LLM-guided framework can generate rules that outperform traditional single priority rules and GP-designed rules obtained under comparable search effort. It also reaches the performance range of the substantially larger GP search with a nominal new-candidate evaluation budget of approximately 2\% of that used by the GP search. Expanding the training set improves performance on all four test sets, while the backbone experiments show that effective rules can be obtained using several different LLMs. For the J30-trained rules, the AllPR multi-rule reference obtains lower AvgDevLB on all five test sets. After expanded training, this ordering reverses on RG300, and the selected LLM rules also outperform AllPR on all three large-project datasets.

The large-project results further highlight the value of reusable priority rules for highly parallel projects, where limited-budget per-instance search can be insufficient. The selected LLM rules outperform AllPR on all three large-project datasets and, on the two highly parallel datasets, also outperform the GA configurations using 1,000 and 5,000 schedule evaluations per instance, while the 50,000-evaluation GA obtains better results. In such cases, a high-quality rule designed offline can provide effective schedules without repeating an extensive search for every project instance.

The full framework performs better than all five ablation variants. This result indicates that the examined search components contribute to the performance of the full framework. The counterfactual analysis further suggests that project-level indicators and progress-dependent switching contributed to the performance of the J30+RG30-trained LLM rules. Algebraically equivalent formulas recur across independent runs and LLM backbones, and rules whose decisions are closer to those of the repeatedly rediscovered rule tend to perform better on the test sets. Behavioral information may therefore help screen candidates, select final and reference rules, and maintain diversity in future searches.

Overall, LLM-guided offline search provides a complementary route to GP-based rule evolution and supervised rule learning. During offline design, the LLM can generate complete executable priority functions without requiring the entire search space to be specified through a fixed grammar. Candidate rules are evaluated directly by the schedules they produce, so the framework does not require labeled schedules or target priority values. The experiments also show that the framework reaches the performance range of the large-budget GP search with far fewer nominal new-candidate evaluations. Once the offline search is complete, the selected rule can be reused on unseen projects as a fixed scoring expression during schedule construction. Deployment therefore requires no LLM calls or iterative search, and the decision logic remains explicit and inspectable.

The evidence in this study is based on benchmark project datasets and makespan minimization. Future research should evaluate the learned rules on empirical project data and in project-planning and rescheduling environments where project information may change during execution and performance criteria extend beyond makespan. The framework could also be adapted to stochastic, multi-mode, and dynamic project-scheduling problems to examine when reusable offline-designed rules remain effective and when online LLM decision support provides additional value. Finally, the observed association between decision behavior and performance motivates studies on using behavioral distance for candidate screening, final-rule selection, reference selection, and diversity management. Such studies could determine whether behavior-based information can improve search efficiency or rule quality when used alongside fitness.

\begingroup
\singlespacing
\section*{Declaration of generative AI and AI-assisted technologies in the manuscript preparation process}

During the preparation of this work, the authors used ChatGPT (OpenAI) to support language editing and manuscript organization. The authors reviewed and edited the content as needed and take full responsibility for the content of the article.

\section*{Funding}

This work was supported in part by IVADO Research Regroupement 9; the Natural Sciences and Engineering Research Council of Canada (NSERC) [grant number 239019]; the industry consortium supporting McGill University's COSMO Stochastic Mine Planning Laboratory; and the Canada Research Chairs Program. The funders had no role in the design or conduct of the study, analysis or interpretation of the results, preparation of the manuscript, or the decision to submit the article.

\bibliographystyle{apalike}
\bibliography{Paper_Ref}
\endgroup

\end{document}